\documentclass{class_for_arxiv}

\usepackage{natbib}
 \NatBibNumeric
 \def\bibfont{\small}%
 \bibpunct[, ]{[}{]}{,}{n}{}{,}%

\usepackage[colorlinks=true,breaklinks=true,bookmarks=true,urlcolor=blue,
     citecolor=blue,linkcolor=blue,bookmarksopen=false,draft=false]{hyperref}

\def\EMAIL#1{\href{mailto:#1}{#1}}

\TheoremsNumberedThrough     

\EquationsNumberedThrough    

\usepackage[english]{babel}

\usepackage{bold-extra}
\usepackage{amssymb}
\usepackage{amsmath}
\DeclareMathAlphabet{\mathpzc}{OT1}{pzc}{m}{it} 
\usepackage[utf8]{inputenc} 
\usepackage[T1]{fontenc}    
\usepackage{hyperref}       
\usepackage{url}            
\usepackage{booktabs}       
\usepackage{nicefrac}       
\usepackage{microtype}      
\usepackage{xcolor}         
\usepackage{dsfont,mathtools}
\usepackage[ruled,vlined]{algorithm2e}
\usepackage{subcaption}
\usepackage{graphicx}
\usepackage{optidef}
\usepackage{enumitem}
\usepackage{diagbox}
\usepackage{tikz-cd}
\usepackage{multirow,arydshln} 

\makeatletter
\newcommand{\mytag}[2]{%
  \text{#1}%
  \@bsphack
  \begingroup
    \@onelevel@sanitize\@currentlabelname
    \edef\@currentlabelname{%
      \expandafter\strip@period\@currentlabelname\relax.\relax\@@@%
    }%
    \protected@write\@auxout{}{%
      \string\newlabel{#2}{%
        {#1}%
        {\thepage}%
        {\@currentlabelname}%
        {\@currentHref}{}%
      }%
    }%
  \endgroup
  \@esphack
}
\makeatother

\graphicspath{{}{figure/}}

\usepackage[textwidth=\linewidth]{todonotes}		

\newcommand\Rmax{r_{\mathrm{max}}}
\newcommand{\pp}{\mathbf{P}}
\newcommand{\rr}{\mathbf{r}}

\newcommand{\bb}{\mathbf{b}}
\newcommand{\bB}{\mathbf{B}}
\newcommand{\bx}{\mathbf{x}}

\newcommand{\calS}{\mathcal{S}}
\newcommand{\calJ}{\mathcal{J}}
\newcommand{\calB}{\mathcal{B}}

\newcommand{\bY}{\mathbf{Y}}

\newcommand{\by}{\mathbf{y}}
\newcommand\xs{\mathpzc{x}}

\newcommand{\R}{\mathbb{R}}

\newcommand{\rel}{\Vrel(\bx,T)}
\newcommand{\Vrel}{V_{\mathrm{rel}}}

\newcommand\ba{\mathbf{a}}
\newcommand\bs{\mathbf{s}}

\newcommand\bX{\mathbf{X}}

\newcommand{\DeltaN}{\Delta^{(N)}}
\newcommand\bXN{\mathbf{X}^{(N)}}
\newcommand\bYN{\mathbf{Y}^{(N)}}

\newcommand\bEN{\mathbf{E}^{(N)}}
\newcommand\XN{X^{(N)}}

\newcommand\YN{Y^{(N)}}

\newcommand\N{\mathbb{N}}
\newcommand\calO{\mathcal{O}}
\newcommand\calE{\mathcal{E}}

\newcommand\calA{\mathcal{A}}
\newcommand\calY{\mathcal{Y}}
\newcommand\calYN{\mathcal{Y}^{(N)}}
\newcommand\calU{\mathcal{U}}

\newcommand{\expect}[1]{\mathbb{E}\left[#1\right]}
\newcommand{\var}[1]{\mathrm{var}\left[#1\right]}
\newcommand{\proba}[1]{\mathbb{P}\left(#1\right)}

\newcommand{\abs}[1]{\left|#1\right|}
\newcommand{\floor}[1]{\left\lfloor#1\right\rfloor}
\newcommand{\norme}[1]{\left\| #1 \right\|_1}
\newcommand{\snorme}[1]{\| #1 \|_1}

\newcommand{\Update}{V^{(N)}_{\mathrm{LP-update}}(\bx,T)}
\newcommand{\UpdateMILP}{V^{(N)}_{\mathrm{MILP-update}}(\bx,T)}
\newcommand{\updates}[2][T-t]{V^{(N)}_{\mathrm{LP-update}}(\bXN(#2),#1)}

\newcommand{\relm}[1]{V_{\mathrm{rel}}(#1)}

\newtheorem{prop}{Proposition}
\newtheorem{lem}{Lemma}
\newtheorem{thm}{Theorem}

\newtheorem{defi}{Definition}

\begin{document}

\RUNAUTHOR{N. GAST, B. GAUJAL, and C. YAN}
%


\TITLE{Reoptimization Nearly Solves Weakly Coupled Markov Decision Processes}

\ARTICLEAUTHORS{
\AUTHOR{Nicolas GAST}
\AFF{Univ. Grenoble Alpes, Inria, \EMAIL{nicolas.gast@inria.fr}}
\AUTHOR{Bruno GAUJAL}
\AFF{Univ. Grenoble Alpes, Inria, \EMAIL{bruno.gaujal@inria.fr}}
\AUTHOR{Chen YAN}
\AFF{Univ. Grenoble Alpes, Inria, \EMAIL{chen.yan@inria.fr}}
}

\ABSTRACT{
We propose a new policy,  called the \emph{LP-update policy}, to solve finite horizon weakly-coupled Markov decision processes. The latter can be seen as multi-constraint multi-action bandits, and generalize the classical restless bandit problems. Our solution is based on re-solving periodically a relaxed version of the original problem, that can be cast as a linear program (LP). When the problem is made of $N$ statistically identical sub-components, we show that the LP-update policy becomes asymptotically optimal at rate $\calO(T^2/\sqrt{N})$. This rate can be improved to $\calO(T/\sqrt{N})$ if the problem satisfies some ergodicity property and  to $\calO(1/N)$ if the problem is \emph{non-degenerate}. The definition of non-degeneracy extends the same notion for restless bandits. By using this property, we also improve the computational efficiency of the LP-update policy. We illustrate the performance of our policy on randomly generated examples, as well as a generalized applicant screening problem, and show that it outperforms existing heuristics.
}



\maketitle

\section{\bfseries\scshape{Introduction}}  \label{sec:introduction}

Markov Decision Processes (MDPs) have proven tremendously useful as models of stochastic sequential planning problems. Dynamic programming is the principal method to address these problems under uncertainty. While MDPs are in theory applicable to many problems, applying classic dynamic programming algorithms to realistic problems is computationally difficult. This has triggered much research into techniques to deal with large state and action spaces. One such technique is decomposition, for which the very large global MDP is decomposed into $N$ loosely dependent sub-processes, each of a size exponentially smaller compared to the global problem. Each local problem can then be solved independently at relatively low cost. If these local solutions can be pieced together effectively, and used to guide the search for a global solution that performs well, then dramatic improvements in the overall solution time can be obtained.

In this paper we study \emph{weakly-coupled} MDPs that fall into this situation. The model originates from sequential stochastic resource allocation problems: A number of tasks must be addressed and an action consists in assigning various resource at every decision epoch to each of these tasks. Weakly-coupled MDPs have been used to model scheduling, queueing, supply chain as well as health care problems  \cite{GOCGUN20122323,hawkins2003langrangian,meuleau1998solving,patrick2008dynamic,salemi2018approximate}. In such a model, the tasks are \emph{additive utility independent}: the utility of achieving any collection of tasks is the sum of rewards associated with each task. Moreover, for a given allocation of resource, each task can be viewed as an independent sub-process whose rewards and transitions are independent of the others. The tasks are linked only through the global resource constraints at each decision epoch. This explains the terminology "weakly-coupled". 

An interesting sub-class of weakly-coupled MDPs, for which many results have been obtained, is the class of multi-armed restless bandits. In such a problem, the bandit (the global MDP) consists of $N$ arms (the sub-MDPs). Two actions are possible for each arm: The passive action and the active action. The decision maker must respect the constraint of activating at most $\alpha N$ arms at each decision epoch. In the seminal paper \cite{whittle-restless}, Whittle proposes a relaxation of the problem that decomposes the $N$-armed problem into $N$ independent sub-problems, one for each arm. Piecing these subproblems together to obtain a feasible solution leads to the famous Whittle index policy, that is known to be extremely efficient. Under several additional technical assumptions, the Whittle index policy is proven to be asymptotically optimal in \citet{WeberWeiss1990}, in the sense that the gap between the performance of the Whittle index policy and the optimal policy converges to zero when $N$ goes to infinity.
In the recent years, many asymptotic optimality results have been obtained for the multi-armed bandit model, under either finite  \cite{Brown2020IndexPA,gast2023linear,zhang2021restless,zhang2022near} or infinite horizons \cite{Ve2016.6,gast2023exponential,zhang2022near}. Some of them assume that arms can have several actions \cite{xiong2021reinforcement,ZayasCabn2017AnAO}.

\paragraph*{Contributions}
We consider a general model of homogenous weakly-coupled MDPs, and we design a policy that becomes optimal as the number of components grows. We provide upper bounds on the sub-optimality gap of the LP-update policy, which is a bound on the rates at which the policy becomes optimal. More precisely, our paper makes two main contributions:
\begin{enumerate}
  \item We generalize the LP-update policy, that first appeared in \citet{gast2023linear} for two-action restless bandits, to the broader class of weakly-coupled MDPs. We call this policy LP-update because it re-solves periodically a new LP to update its future control decisions.  We show that the sub-optimality gap of the LP-update policy is at most $\calO(T^2/\sqrt{N})$, where $N$ is the number of sub-components of the weakly-coupled MDP and $T$ is the time horizon. Under the condition that the ergodic coefficient of each homogenous sub-component is strictly smaller than one, we show that the sub-optimality gap can be refined to $\calO(T/\sqrt{N})$. These rates improve upon the bound $\calO(T^3/\sqrt{N})$ in \cite{brown2023fluid}, under the more general case of heterogeneous sub-components.
      
  \item We introduce the notion of a \emph{non-degenerate} weakly-coupled MDP using a rank condition inspired from the linear independence constraint qualification in optimization (\citet{bazaraa2013nonlinear}), and show that for such a problem, the sub-optimality gap of the LP-update policy is $\calO(1/N)$. If the problem satisfies an additional perfect rounding condition, then the rate can be further improved to $e^{-\Omega(N)}$. The non-degeneracy is a key notion to achieve faster convergence rates. There already exists definitions of non-degeneracy for particular sub-classes of weakly-coupled MDPs: The ones presented in \cite{gast2023linear,zhang2021restless,brown2023fluid} focus on particular sub-cases, and are equivalent to our definition when restricted to those. The one given in \citet{zhang2022near1} for general weakly-coupled MDPs is more restrictive than ours, in the sense that the class of non-degenerate problems according to \citet{zhang2022near1} is a strict subset compared to using our definition. This difference is important in practice: on randomly generated examples, about half of our non-degenerate examples are not captured by the definition of \cite{zhang2022near1}.
\end{enumerate}
In addition, we show that all the convergence rates claimed previously are tight, in the sense that for each case, we are able to exhibit an example whose sub-optimality gap does not decrease faster than our upper bound.  This shows that we have identified the essential features of the problem for achieving a fast convergence rate.
We also present a few numerical results that illustrate the benefit of using our LP-update policy compared with the occupation measure policy presented in \cite{xiong2021reinforcement,ZayasCabn2017AnAO}. We illustrate this on randomly generated examples, as well as on a case study that is a generalization of the applicant screening problem studied in \citet{Brown2020IndexPA} and \citet{gast2023linear}. 

\paragraph*{Outline}
The rest of the paper is organized as follows. We discuss related work in Section~\ref{sec:related-work}. We introduce the weakly-coupled MDPs model in Section \ref{sec:model-description}. The LP-update policy for weakly-coupled MDPs with full updates is given in Section \ref{sec:LP-update-in-general}. The improved version of the policy with selective updates is given in Section \ref{sec:non-degenerate}, together with the definition of non-degeneracy. The case study on the applicant screening problem is given in Section \ref{sec:general-applicant-screening}. We conclude the paper and discuss some future research directions in Section \ref{sec:conclusion-and-future-work}. The proofs of the technical results are collected in Section \ref{sec:proof-of-main-theorem}. Additional discussions are collected in the other appendices.

\section{Related work}
\label{sec:related-work}

There are several branches of works in the existing literature that are related to our results in this paper. The first branch is concentrated on the study of weakly-coupled MDPs under finite horizon, e.g. \citet{10.1287/opre.1070.0445}, \citet{ASTARAKY2015309}, \citet{dolgov04resourceMDP}, \citet{5679024}, \citet{GOCGUN20122323}, \citet{meuleau1998solving}, \citet{patrick2008dynamic}, as well as the two PhD thesis \citet{hawkins2003langrangian} and \citet{salemi2018approximate}. In these works the authors consider the situation when the sub-MDPs are not statistically identical. The authors of these papers use a Lagrange decomposition technique to solve approximately the above problem, which works as follows: By relaxing the constraints with state independent multipliers, the $N$-dimensional optimization problem decouples into $N$ one-dimensional subproblems, each involves the corresponding sub-MDP alone; the key is in choosing the value of multipliers and in transforming the optimal controls of the decoupled problems into a feasible policy to the original problem (as has been done for the simpler multi-armed bandit model). The idea of decomposition also appears in other domains, for instance in network congestion control as in \citet{palomar2006tutorial}. More broadly, in a series of papers \cite{carpentier2018stochastic, carpentier2020mixed, pacaud2021distributed, pacaud2022optimization} as well as in the textbook \citet{carpentier2017decomposition}, the authors study a very general large-scale stochastic optimization problem with many distributed sub-components, using a guiding principle called "decomposition-coordination". 

Our results are also close to  the branch of researches on asymptotic optimality results on multi-armed bandits under finite horizon as in \citet{Brown2020IndexPA}, \citet{zhang2021restless}, \citet{gast2023linear}, \citet{ZayasCabn2017AnAO} and \citet{xiong2021reinforcement}. The finite horizon two-action single-constraint multi-armed bandit is the subject of the first three papers, the more general multi-action single-constraint multi-armed bandit is the subject of the two latter papers. Except in \citet{gast2023linear}, all the policies considered in these references are one-pass policies that involve solving a linear program only once at the very beginning. The non-degenerate property on two-action single-constraint multi-armed bandit is proposed independently in both \citet{zhang2021restless} and \citet{gast2023linear}, together with the $\calO(1/N)$ rate on non-degenerate models proven therein. 

There is another related branch of works that consider the problem under infinite horizon. In \citet{whittle-restless}, \citet{WeberWeiss1990}, \citet{Ve2016.6}, \citet{gast2023exponential} and \citet{zhang2022near} the two-action bandit has been studied.
It is shown in \citet{WeberWeiss1990} that, for indexable model, Whittle index policy is asymptotically optimal under an additional global attractor assumption. Later in \citet{Ve2016.6} the indexability assumption has been removed. In \citet{gast2023exponential} the optimal convergence rate of these asymptotic results are proven. In \citet{zhang2022near}, by using the discount criterion, instead of the time-average criterion considered in other works, the authors construct an asymptotically optimal policy that does not require the indexability nor the global attractor assumption. This asymptotic result relies on a truncation over a finite horizon $T$, that depends on the discount factor. When studying the problem under finite horizon as we do here, none of the assumptions of indexability nor the global attractor property are needed for the asymptotic optimality results. On the other hand, the policies become time-dependent and the focus shifts to the transient behavior. The notion of index is generalized to the multi-action single-constraint case in \citet{hodge2015asymptotic} and \citet{math10142497}. It is not clear if such a notion can be generalized to the case of multiple constraints, since the index of a sub-process can only be defined with respect to a single constraint, and it is not obvious to incorporate multiple constraints into one real-valued meaningful index. 
The infinite-horizon two-action bandit problem is also explored in the recent study by \cite{hong2023restless}, wherein the authors present an algorithm named \emph{Follow-the-Virtual-Advice} (FTVA). Assuming a mild condition termed synchronization by the authors, the suboptimality gap of FTVA is shown to be $\Rmax\tau^{\mathrm{sync}}/\sqrt{N}$, where $\tau^{\mathrm{sync}}$ represents a coupling time. This synchronization condition bears resemblance to our ergodicity condition $\lambda^{\tau}>0$ in Theorem~\ref{thm:LP-update-general-refined}. Despite our analysis focusing on the finite horizon with multiple constraints, as opposed to their single constraint infinite horizon framework, direct comparison of the two performance bounds is challenging. On one side, our bound incurs an additional factor of $\sqrt{d}$; on the other, their bound is influenced by a "passive" coupling time of independent trajectories, while ours relies on an "active" mixing time typically much shorter than the coupling time for large $d$.

Our proposed heuristic for weakly-coupled MDPs, that we call the LP-update policy, generalizes the one for two-action single-constraint bandits proposed in \citet{gast2023linear}, with one difference: in the current paper, we improve its computational complexity by re-solving a LP only when necessary. We also obtain a faster rate of convergence than the $\calO(1/\sqrt{N})$ proven in \cite{gast2023linear}. The LP-update policy combines ideas from the certainty equivalent control in dynamic programming (Chapter 6 of \citet{bertsekas2012dynamic}), with the model predictive control in optimal control theory (\citet{rakovic2018handbook}). The same technique has been applied in other stochastic optimization problems. For example, the discrete-review policy for multi-class queueing network in \citet{maglaras2000discrete}, or the re-solving heuristic for the network revenue management in \citet{jasin2012re}.  In this paper we use a similar idea as the ones used in these papers that consists of four ingredients: (a) find an appropriate deterministic approximation to the stochastic problem (see Section \ref{subsec:relaxed-problem}); (b) use the optimal solution to the deterministic problem to derive an implementable solution for the original problem (see Section \ref{subsec:motivation-for-updates} and Section \ref{subsec:rounding-problem-discussion}); (c) decide the moments for re-optimization (see Section \ref{subsec:improved-LP-update}); (d) analyse the performance of the resulting policy (see Theorems \ref{thm:LP-update-general}, \ref{thm:LP-update-general-refined}, \ref{thm:LP-update-rate-nondegenerate} and \ref{thm:LP-update-rate-expo}).

Concurrently with the development of this study, an independent research project led by \citet{brown2023fluid} explores a related issue. Their work primarily delves into heterogeneous weakly-coupled MDPs, contrasting with our investigation of homogeneous sub-MDPs. The convergence rate identified by \citet{brown2023fluid} for their "reoptimized feasible fluid policy" in a broader context is $\calO(T^3/\sqrt{N})$, where $T$ is the finite horizon and $N$ is the number of sub-components. In comparison, our approach yields a convergence rate of $\calO(T^2/\sqrt{N})$, and potentially $\calO(T/\sqrt{N})$ when the ergodic coefficient of the sub-components is less than one. This improved convergence rate in our model relies on the premise of homogeneity across all arms. Furthermore, our bounds exhibit simple (essentially linear) dependence on additional problem parameters, such as the cardinality of the state and action space, as articulated in Theorem \ref{thm:LP-update-general}. This is distinct from bounds of the theorems of Section~4 of \cite{brown2023fluid} that have more complex dependence on the problems parameters.


\section{\bfseries\scshape{Model Description}} \label{sec:model-description}

\subsection{Weakly-coupled MDPs}

We consider a finite-horizon discrete-time weakly-coupled MDP composed of $N$ statistically identical sub-MDPs (or \emph{components}, or \emph{arms} in the bandit literature), indexed by $n\in\{1\dots N\}$. The finite state space of each sub-MDP is the set $\calS := \{1,2,\dots,d \}$, and its finite action space\footnote{For notational simplicity, we assume that the sub-action space of each sub-MDP $n$ is independent of its state $S_n$. This is without loss of generality as one may forbid some actions by giving them an extremely high cost.} is $\calA:=\{0,1,\dots,A\}$. The state space of the weakly-coupled MDP is therefore $\calS^N$ and the action space is a subset of $\calA^N$.  There are $J$ types of resources, and the decision maker is allowed to use up to $N \cdot b_j>0$ resource of type $j$ at each decision epoch.  We assume that taking the action $a_n$ for component $n$ being in state $s_n$ uses $D_j(s_n,a_n)\ge0$ of resource $j$, and that the action $0$ consumes no resource: $D_j(s_n,0)=0$ for all $s_n\in\calS$.  Hence, the set of \emph{feasible} actions in state $\bs = (s_1,s_2, \dots, s_N)$ is the set of $\ba \in \calA^N$ such that for all $j\in\{1\dots J\}$: $\sum_{n=1}^{N} D_j(s_n,a_n) \le N \cdot b_j$.

The sub-processes are weakly-coupled, in the sense that the $N$ sub-MDPs are only linked through the $J$ resource constraints, i.e. for a given feasible action $\ba$, the system transitions from a state $\bs$ to state $\bs' = (s'_1,s'_2, \dots, s'_N)$ with probability
\begin{equation}\label{eq:Markov-transition-each-arm}
  p(\bs' \mid \bs,\ba) = \prod_{n=1}^{N} p(s'_n \mid s_n,a_n) = \prod_{n=1}^{N} P^{a_n}_{s_n, s'_n},
\end{equation}
where for each action $a$, the matrix $\pp^a$ is a probability transition matrix of dimension $d \times d$.

Upon choosing an action $\ba$ in state $\bs$, the decision maker receives an immediate reward $\sum_{n=1}^{N} r_{s_n}^{a_n}$, where $r_{s_n}^{a_n} \in \mathbb{R}$ depends on the state $s_n$ and action $a_n$, and the rewards among components are additive. We denote by $\Rmax:=\max_{s,a}\abs{r^{s,a}}<\infty$. It is finite as the state and action spaces are finite. 

The aforementioned model makes a number of assumptions, that are classical in the literature (e.g. \cite{xiong2021reinforcement,ZayasCabn2017AnAO}):
\begin{itemize}
  \item We assume that all terms in $D(s,a)$ and $\bb$ are non-negative numbers, and that $D(s,0)= \mathbf{0}$. This is a natural assumption under the resource allocation context in which $a=0$ corresponds to a \emph{passive} action that consumes no resource. The later also implies that our resource constraint problem has at least a feasible solution by always choosing the passive action. .
  \item We assume that each sub-MDP is statistically identical. This is needed under our scaling of the arm population $N$. However, this assumption can be relaxed, in the sense that we can incorporate the case where there is a finite number of types of sub-MDPs, by making direct sum of the state spaces of each type into a single larger state space.
  \item In our formulation of the problem, the rewards, constraints and transition probabilities do not dependent on time.  This choice is only to lighten the notations: As we consider finite-horizon problem, all the results apply to the case of time-dependent parameters---it would suffice to add a dependence on $t$ on all parameters.
\end{itemize}
Finally, we remark that this model includes the classical restless multi-armed bandit that corresponds to the case $\calA=\{0,1\}$ with $J=1$ constraint and $D_1(s,a)=a$. This simpler optimization problem is already PSPACE-hard (see  \citet{Papadimitriou99thecomplexity}). In the rest of the paper, we focus on developing approximate solutions whose performance are provably close to optimal.

\subsection{Symmetry simplification and population representation}

Since we assume that the $N$ sub-components are statistically identical, the problem can be reformulated by keeping track of the number (population) of sub-MDPs in each of the $d$ states, as well as the number of sub-MDPs for each action taken in that state. For the purposes of our subsequent analyses, it proves beneficial to normalize each quantity by dividing it by $N$.

Hence, we denote by $\bXN(t) = (\XN_s(t))_{s\in \calS}$  the vector of proportions of sub-MDPs that are in each state $s\in \calS$ at decision epoch $t$. We denote by $\DeltaN$ the possible values for $\bXN$, which is the set of vectors $\bx\in\R^d$ such that $x_s\ge0$, $\sum_{s\in\calS}x_s=1$ and $N \cdot x_s$ is an integer for all $s\in\calS$. It is a finite collection of points in the simplex $\Delta^d$. Although this representation reduces the state space size of the MDP considerably from $d^N$ to $\binom{N+d-1}{d-1}$, it is still intractable when $N$ and $d$ are large.

Let $\calU=\calS\times\calA$ be the set of state-action pairs. Upon observing $\bXN(t)$, the action taken by the decision maker is represented by  $\by\in \R^u$ where $y_{s,a}$ is the proportion of sub-MDPs that are in state $s$ and for which action $a$ is taken.  For a given $\bx\in\DeltaN$, we denote by $\calYN(\bx)$ the set of feasible actions in state $\bx$. It is the set of possible $\by\in\calY(\bx)$ such that for all $s,a$, $N \cdot y_{s,a}$ is an integer, where $\calY(\bx)$ is:
\begin{align}
  \label{eq:feasible_action}
  \calY(\bx) &:= \left\{ \by\ge0 \text{ such that } \forall s\in\calS: \sum_{a\in\calA} y_{s,a}=x_s;  D \by \le \bb \right\}.
\end{align}
In the above definition, the first set of equality constraints is due to the mass conservation, and the second represents the budget constraints. They are written in matrix form by viewing $D_{j}(s,a)$ as a matrix whose lines are indexed by $j\in\{1\dots J\}$ and whose columns are indexed by $(s,a)$. Similarly, the vector $\bb$ is a column vector indexed by $j\in\{1\dots J\}$.

\subsection{Optimal control formulation}


The decision maker's goal is to maximize the total expected reward over a finite horizon $T$, given an initial state vector $\bx \in \DeltaN$ of the system. We denote by $V^{(N)}_{\mathrm{opt}}(\bx,T)$ the maximal expected gain per arm that can be obtained by the decision maker.

This optimization problem can be written as follows:
\begin{maxi!}|s|{\bY}{\mathbb{E} \Big[ \sum_{t=0}^{T-1} \sum_{(s,a) \in \calU} r^{a}_{s} \YN_{s,a}(t) \Big] }{\label{eq:prob-formulated}}{V^{(N)}_{\mathrm{opt}} (\bx,T) =}
  \addConstraint{\XN_s(0) = x_s \ \ \forall s \label{eq:prob-formulated_M0} }
  \addConstraint{\text{$\bXN(t+1)$ follows the Markov transitions \eqref{eq:Markov-transition-each-arm} given  $\bYN(t)$ }  \label{eq:prob-formulated_Mt+1}}{}
  \addConstraint{\bYN(t) \in \calYN(\bXN(t)) \ \ \forall t \label{eq:feasible-action-stage-t} }{}
\end{maxi!}

\subsection*{Notation and Terminology Convention}
A super-script $N$ indicates that the problem has population $N$, with $N$ being a positive integer number. This implies  that variables with super-script $(N)$ (for example $\XN_s$) are integral multiples of $1/N$ (i.e. $N \cdot \XN_s \in \N$). When we talk about asymptotic results, this variable $N$ is understood as a scaling parameter that tends to infinity. By "asymptotic optimality", implicitly we are always comparing with the LP relaxation bound, which is an upper bound on the value of the optimal policy. We do so since the later is in general not easily computable. A bold letter (e.g. $\by$, $\bx$) denotes a vector whereas a normal letter (e.g. $y_{s,a}(t)$, $x_s(t)$) denotes a scalar. Unless otherwise specified, vectors are column vectors. We use super-script ${}^*$ to denote optimal quantities (e.g. $\bx^*$, $y^*$). The function $\mathbf{1}_{E}$ is a random variable that equals $1$ if the event $E$ occurs and $0$ otherwise. The vector $\rr^\top$ denotes the transpose of the vector $\rr$. For two real valued functions $f$ and $g$, we write $f(x) = \calO(g(x))$ if there exists $C>0$ and $x_0$ real such that $\abs{f(x)} \le C g(x)$ for all $x \ge x_0$, and we write $f(x) = \Omega(g(x))$ if $\limsup_{x \rightarrow \infty} \abs{f(x) / g(x)} > 0$.

\section{\bfseries\scshape{The LP-update Policy for Weakly Coupled MDPs}} \label{sec:LP-update-in-general}

In this section, we first introduce the linear program after relaxing the resource constraints in Section \ref{subsec:relaxed-problem}, which is the starting point for all policies we consider afterwards. We then define the LP-update policy with full updates, and present the $\calO(1/\sqrt{N})$ performance guarantee that holds in all cases in Section~\ref{sec:LP-update}.

\subsection{The relaxed problem as a linear program} \label{subsec:relaxed-problem}

The main difficulty of the optimization problem \eqref{eq:prob-formulated} is that the constraint $D\bYN(t)\le \bb$ couples all sub-MDPs. To overcome this difficulty, a by-now classical approach \cite{Brown2020IndexPA,gast2023linear,xiong2021reinforcement,ZayasCabn2017AnAO,zhang2021restless} is to relax these constraints and consider a problem where they have to be satisfied only in expectation: $D\expect{\bYN(t)}\le \bb$.  This leads us to write a relaxed optimization problem in terms of the variables $\by(t)=\expect{\bYN(t)}$.  Indeed, equation \eqref{eq:Markov-transition-each-arm} implies that the expectation of $\bXN(t+1)$ given $\bYN(t)$ can be rewritten as a linear map $\phi$ as follows:
\begin{equation}\label{eq:Markov-transition-expected}
  \expect{\XN_{s}(t+1)\mid \bYN(t)=\by} = (\phi(\by))_s := \sum_{(s',a) \in \calU} y_{s',a} P^a_{s',s}.
\end{equation}
This shows that the relaxed optimization problem is the following linear program with decision variables $\by(t)=\expect{\bYN(t)}$:
\begin{maxi!}|s|{\by\ge0}{\sum_{t=0}^{T-1} \sum_{(s,a) \in \calU} r_s^a y_{s,a} (t)}{\label{eq:relaxed_problem_general}}{V_{\mathrm{rel}}(\bx,T)=}
  \addConstraint{\sum_{a \in \calA} y_{s,a}(0) = x_s}{}{\forall s \label{eq:init_general}}
  \addConstraint{\sum_{a \in \calA} y_{s,a}(t+1) = \big( \phi(\by(t)) \big)_s \label{eq:markov_relaxed_general}\qquad}{}{\forall s,t}
  \addConstraint{D \by(t) \le \bb \label{eq:relaxed-constraints}}{}{\forall t},
\end{maxi!}
 In the above formulation, the constraint \eqref{eq:init_general} corresponds to the condition on the initial state. The constraint \eqref{eq:markov_relaxed_general} corresponds to the time-evolution \eqref{eq:Markov-transition-each-arm}, plus the fact that $x_s(t+1)=\sum_a y_{s,a}(t+1)$. Lastly, the constraint \eqref{eq:relaxed-constraints} corresponds to $D\expect{\bYN(t)}\le \bb$, which is the relaxed version of the one in \eqref{eq:feasible_action}. Equivalently it can be written as $\by(t)\in\calY(\bx(t))$, where $x_s(t)=\sum_{a\in\calA}y_{s,a}(t)$. For later reference, we call the linear program \eqref{eq:relaxed_problem_general} being parameterized by $(\bx,T)$, where $\bx$ is the initial state vector and $T$ is the horizon. 

By the assumptions that $D(s,0)= \mathbf{0}$ and $D(s,a), \bb \ge \mathbf{0}$, the linear program \eqref{eq:relaxed_problem_general} is feasible (\emph{e.g.} it suffices to always choose the passive action $a=0$). In the following, we denote by $\by^*$ one of its optimal solution, and by $\bx^*$ the sequence of vectors $\bx^*(t)$ such that $x^*_s(t)=\sum_{a\in\calA}y^*_{s,a}(t)$. It is the optimal state vector on the relaxed problem at time-step $t$.

\subsection{The LP-update policy with full updates} \label{subsec:motivation-for-updates} \label{sec:LP-update}

To construct a policy for the system of size $N$, the LP solution suggests to use $\bYN(t)=\by^*(t)$ at decision epoch $t$. Yet, this is in general not possible because of random fluctuations: Indeed, it is likely that $\bXN(t)\ne\bx^*(t)$, which implies that in general, $\by^*(t)$ is not feasible for $\bXN(t)$, that is
\begin{align}
  \label{eq:non-feasible}
  \text{ In general: } \by^*(t)\not\in\calYN(\bXN(t)).
\end{align}
The classical way to solve this problem in the literature is to construct a sequence of decision rules $\pi_t:\Delta^d\to\Delta^{d(A+1)}$ such that $\pi_t(\bx)\in\calYN(\bx)$ and $\pi_t(\bx^*(t))=\by^*(t)$. This is what has been used to build the randomized activation control policy in \citet{ZayasCabn2017AnAO}, the fluid-priority policies in \citet{zhang2021restless}, the LP-index policy in \citet{gast2023linear}, the occupancy-measured-reward index policy in \citet{xiong2021reinforcement}, to name a few.  In particular, it is shown in \cite{gast2023linear} that any such  policy is $\calO(1/\sqrt{N})$-optimal if all the decision rules $\pi_t(\cdot)$ are locally Lipschitz-continuous functions.

In this paper, we adopt another approach, that we call the \emph{LP-update} policy, which is a generalization of the LP-update policy of \citet{gast2023linear} introduced for two-action bandits. The policy of \cite{gast2023linear} is as follows: at each decision epoch, we solve a new LP starting from $\bXN(t)$ with horizon $T-t$. This guarantees that the newly computed $\bY(t)$, that corresponds to the solution $\by^*(0)$ of the LP, is in $\calY(\bXN(t))$, by the mass conservation \eqref{eq:init_general}. However, this control is not necessarily feasible for the system of size $N$, since $N \cdot Y_{s,a}(t)$ is not necessarily an integer. A simple\footnote{We will discuss more sophisticated rounding procedures in Section~\ref{subsec:rounding-problem-discussion}.} way to obtain a feasible solution is to use a rounding procedure: 
\begin{align}
  \label{eq:rounding}
  \YN_{s,a}(t) =  \left\{\begin{array}{ll}
    N^{-1}\floor{N \cdot Y_{s,a}(t)} & \text{ if $a\ne0$.}\\
    \XN_{s}(t) - \sum_{a>0} \YN_{s,a}(t)  &\text{ if $a=0$.}
  \end{array}\right.
\end{align}
which satisfy the resource constraints because the resource consumptions of action $0$ are all $0$.

This leads to our first LP-update algorithm, that is summarized in Algorithm~\ref{algo:LP-update}.

\begin{algorithm}[htbp]
  \SetAlgoLined
  \SetKwInput{KwInput}{Input}
  \KwInput{Initial configuration vector $\bXN(0)$ over  time span  $[0,T]$.}
  \For{$t = 0,1,2,\dots,T-1$}{
    Solve LP \eqref{eq:relaxed_problem_general} with parameters $(\bXN(t),T-t)$. Output is $\by^*$ over $[0,T-t]$ \;
    Set $\YN_{s,a}(t) := N^{-1}\floor{Ny^*_{s,a}(0)}$ for $a\not= 0$ and  $\YN_{s,0}(t) := \XN_{s}(t) - \sum_{a\not = 0} \YN_{s,a}(t)$ \;
    Use actions $\YN_{s,a}(t)$ over all sub-MDPs to advance to the next time-step\;
  }
\caption{The LP-update policy for weakly-coupled MDPs (full updates).}
  \label{algo:LP-update}
\end{algorithm}

\subsection{Performance guarantees: general case}

The forthcoming Theorem~\ref{thm:LP-update-general} demonstrates that the LP-update policy, as outlined in Algorithm \ref{algo:LP-update}, achieves an $\calO(T^2/\sqrt{N})$-optimal performance. This optimality further improves to $\calO(T/\sqrt{N})$ when the probability transition matrices meet a certain ergodic condition. The proof of this theorem is anchored on two pivotal elements. Firstly, it leverages the characteristics of the LP-update algorithm to establish a connection between the sub-optimality gap and the Lipschitz coefficient of the value function $\bx \mapsto \Vrel(\bx,t)$. Secondly, it derives an explicit bound on this Lipschitz coefficient. 

It is important to note that the first phase of this proof is inapplicable to LP-based policies that do not involve resolving the LP at every step. This distinction is crucial, as it is the resolution of the LP at each step in our approach that enables us to achieve a polynomial dependence on the time horizon $T$ for the optimality gap.

\begin{thm} \label{thm:LP-update-general} \
  Denote by $\Update$ the value of the expected performance obtained by applying Algorithm~\ref{algo:LP-update} 
  and by $\rel$ the value of the linear program \eqref{eq:relaxed_problem_general}.
 \begin{enumerate}
   \item[(i)] For any weakly-coupled MDP with statistically identical arms and all $N$, we have:
   \begin{equation*}
     \abs{\Update - \rel} \le \frac{\Rmax\sqrt{d}}{4}\frac{T^2}{\sqrt{N}} +\frac{d|\calA|T\Rmax}{N}.
   \end{equation*}
    \item[(ii)] For any weakly-coupled MDP with statistically identical arms such that\footnote{The value of $\gamma$ is very similar to what is call the ergodic coefficient of a MDP in \cite{Puterman:1994:MDP:528623}. The ergodic coefficient contains an additional minimization on action and does not impose $i \neq j$. It is equal to $\min_{i,j, a, b} \sum_k \min( P^{b}_{ik}, P^a_{jk})$ and is smaller or equal to $\gamma$.} $\gamma = \min_{i \neq j, a} \sum_k \min( P^0_{ik}, P^a_{jk})>0$, we have for all $N$:
    \begin{equation*}
     \abs{\Update - \rel} \le \frac{\Rmax\sqrt{d}}{2\gamma}\frac{T}{\sqrt{N}} +\frac{d|\calA|T\Rmax}{N}.
   \end{equation*}
   \item[(iii)] There exists a weakly-coupled MDP with statistically identical arms and a constant $C'>0$ such that for all $N$:
   \begin{equation*}
     \abs{\Update - \rel} \ge \frac{C'}{\sqrt{N}}.
   \end{equation*}
 \end{enumerate}
\end{thm}

\begin{proof}{Proof.}
  The proof of the lower bound (iii) is done in Section~\ref{subsec:lowerbound}. Below, we prove Items (i) and (ii).

  We start our analysis with Algorithm~\ref{algo:LP-update} that executes an update at every time step. Specifically, at each decision epoch $t$, the LP-update algorithm selects a vector $\bY(t)$ deemed optimal for the LP defined in \eqref{eq:relaxed_problem_general}, using the parameters $(\bXN(t), T-t)$. Following this, a decision vector $\bYN(t)$ is chosen from the set $\calYN \left(\bXN(t)\right)$. This vector $\bYN(t)$ is designed to be close to $\bY(t)$, with a deviation bounded by $\norme{\bYN(t)-\bY(t)}\le c/N$, as stipulated in our construction in \eqref{eq:rounding}. Upon implementing this decision, the controller receives a reward equal to $\rr^\top \bYN(t)$. Subsequently, the system progresses to the next state $\bX(t+1)$. This process translates into
  \begin{align}
    \label{eq:proof_update}
    \updates[T-t]{t} &= \expect{\rr^\top \bYN(t) + \updates[T-t-1]{t+1}}.
  \end{align}
 On the other hand,  Bellman's principle of optimality says 
  \begin{align}
    \label{eq:proof_rel}
    \relm{\bXN(t),T-t} &= \rr^\top \bY(t) + \relm{\phi(\bY(t)),T-t-1}.
  \end{align}
   Denote by $Z(t) := \expect{\updates[T-t]{t}- \relm{\bXN(t),T-t}}$ for $0 \le t \le T+1$, where we interpret $Z(T) = Z(T+1) = 0$. Combining \eqref{eq:proof_update} and \eqref{eq:proof_rel} shows that $Z(t)-Z(t+1)$ is equal to
  \begin{align}
    \label{eq:Z(t)}
    \underbrace{\rr^\top \expect{\bYN(t)-\bY(t)}}_{\mytag{Term A}{term:A}} +  \underbrace{\expect{\relm{\bXN(t+1),T-t-1} -  \relm{\phi(\bY(t)),T-t-1}}}_{\mytag{Term B}{term:B}}.
  \end{align}
  By construction, \ref{term:A} is smaller than $\Rmax d |\calA|/N$. Moreover, we shall show in Lemma~\ref{lem:Lipschitz-lemma} that the function $\bx \mapsto V_{\mathrm{rel}} ( \bx , t)$ is $L_t$-Lipschitz-continuous, with $L_t\le\min(t/2, 1/(2\gamma))$. Hence, the absolute value of \ref{term:B} is such that
  \begin{align*}
    \abs{\eqref{term:B}} \le \min\Big(\frac{T-t-1}{2}, \frac1{2 \gamma}\Big) \expect{\norme{\bXN(t+1) - \phi(\bY(t))}}.
  \end{align*}
  By Lemma~\ref{lem:Markovian-transition-analysis}, the inequality
  $\expect{\norme{\bXN(t+1) - \phi(\bY(t))}}\le\sqrt{d}/\sqrt{N}$ holds. The result then follows by summing over $t$ from $t=0$ to $t=T$ for the differences $Z(t)-Z(t+1)$, depending on whether $\gamma > 0$ or not.
  \qed
\end{proof}

Theorem~\ref{thm:LP-update-general}(ii) shows that the sub-optimality gap of the LP-update policy is at most linear in $T$ when $\gamma>0$. The quantity $\gamma$ corresponds to the minimum probability that two components starting in two different states couple after one time-step, where one component takes the action $0$ and the other takes an arbitrary action. In fact, it is possible to generalize this notion by considering the probability that two components couple after more than one time step. More precisely, for any $\tau\in\{1,2\dots\}$, we define $\lambda^\tau$ as:
\begin{align} \label{eq:lambda_tau}
  \lambda^\tau := \min_{i \neq j, a_1\dots a_\tau} \sum_k \min( ( \underbrace{P^0\dots P^0}_{\tau \text{ times}} )_{ik}, (P^{a_1}\dots P^{a_\tau})_{jk})
\end{align}
Note that a sufficient condition for the existence of a $\tau$ such that $\lambda^{\tau}>0$ is that the Markov chain corresponding to $P^0$ is irreducible and aperiodic. Indeed, by classical Markov chain theory, ergodicity and aperiodicity of $P^0$ implies that there exists a time $\tau\ge1$ such the matrix $(P^0)^\tau$ has all of its coefficient strictly positive. This implies that $\lambda^\tau>0$ for this value of $\tau$. We obtain the following theorem, that is a generalization of Theorem~\ref{thm:LP-update-general}(ii). 
\begin{thm}  \label{thm:LP-update-general-refined}
  For any weakly-coupled MDP with statistically identical arms such that there exists an integer $\tau \ge 1$ such that  $\lambda^\tau>0$, we have:
   \begin{equation*}
    \abs{\Update - \rel} \le \frac{\tau\Rmax\sqrt{d}}{2\lambda^\tau}\frac{T}{\sqrt{N}} +\frac{d|\calA|\Rmax T}{N}.
  \end{equation*}
\end{thm}
\begin{proof}{Proof.}
  The proof of this result is identical to the one of Theorem~\ref{thm:LP-update-general}(ii), but we replace the use of Lemma~\ref{lem:Lipschitz-lemma} by the one of Lemma~\ref{lem:Lipschitz-lemma-tau} (these lemma are given in Section \ref{subsec:lipschitz-const}).  
\end{proof}


Theorem~\ref{thm:LP-update-general}(i) shows that our algorithm offers a performance guarantee of $\calO(1/\sqrt{N})$. As established in point (iii) of the Theorem, this rate is optimal for general models and cannot be further enhanced within that context. However, in Section \ref{sec:non-degenerate}, we will introduce the concept of \emph{non-degenerate} problems. This concept builds upon and extends the ideas presented by \citet{gast2023linear,zhang2021restless}. For these specific non-degenerate problems, we demonstrate that the sub-optimality gap can be improved to $\calO(1/N)$. Moreover, we will explore how leveraging the non-degenerate property of problems not only refines our theoretical understanding but also yields computational advantages.
\subsection{Resolving LP brings performance benefits}

Algorithm~\ref{algo:LP-update} solves a new LP at each decision epoch. This can be computationally intensive, making its time complexity high compared with the methods proposed by \cite{Brown2020IndexPA,gast2023linear,xiong2021reinforcement,ZayasCabn2017AnAO,zhang2021restless}. In this subsection, we show that this time cost is compensated by its performance gain. For that, we focus on a numerical comparison between our LP-update policy and the \emph{occupation measure} policy, introduced  in \cite{ZayasCabn2017AnAO,xiong2021reinforcement}. This algorithm solves an LP only once and then bases all subsequent decisions on this initial solution. The occupation measure policy has an asymptotic sub-optimality gap in $\calO(1/\sqrt{N})$, as established in \cite{ZayasCabn2017AnAO}. 

However, a closer look at the proof of the occupation measure's optimality shows that its hidden constant in the $\calO(\cdot)$ notation grows exponentially with the time horizon $T$, due to estimations similar to those in Gronwall's lemma. In contrast, the sub-optimality gap of our LP-update algorithm has  a more modest growth rate of $\calO(T^2/\sqrt{N})$ or $\calO(T/\sqrt{N})$.  

\begin{figure}[ht]
  \centering
  \begin{tabular}{cc}
    \includegraphics[width=0.45\linewidth]{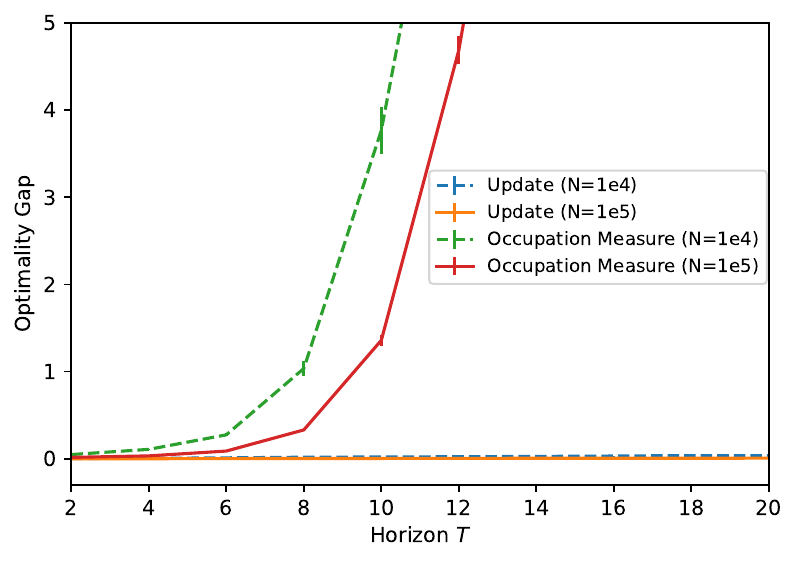}
    &\includegraphics[width=0.47\linewidth]{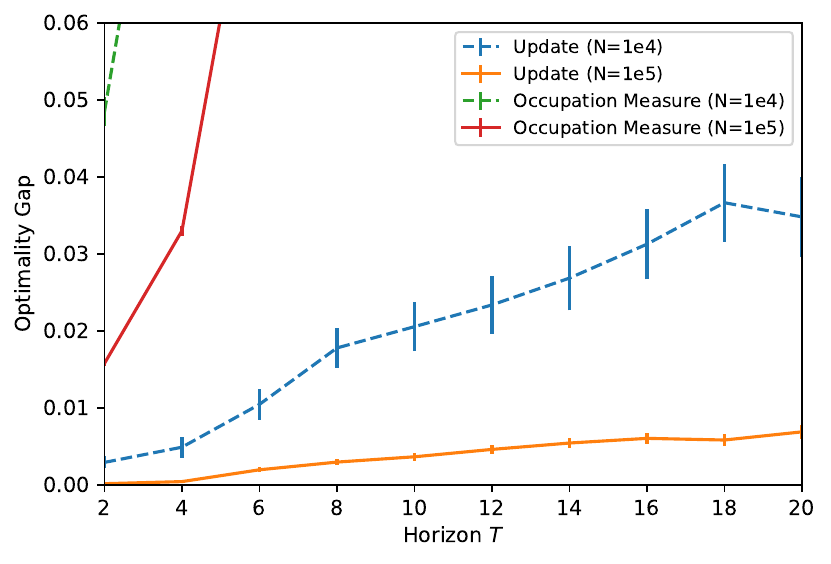}\\
    (a) LP-update vs occupation measure
    &(b) Zoom of the left panel.
  \end{tabular}
  \caption{Sub-optimality gap as a function of the time-horizon $T$. The right panel is a zoom-in of the left panel.}
  \label{fig:growth-rate}
\end{figure}

In Figure~\ref{fig:growth-rate}, we illustrate the sub-optimality gap of both policies as a function of the time horizon $T$ for a model generated at random. This model assumes a large number of identical arms, with $N$ being either $10^4$ or $10^5$. Our observations indicate that the sub-optimality gap for the occupation measure policy appears to increase exponentially with $T$. In contrast, the sub-optimality gap of the LP-update algorithm displays a linear growth relative to $T$. This finding underscores that solving a new LP at each time step can offer significant performance advantages over an algorithm that resolves the LP only once at the initial step. 

In the upcoming section, we will explore how the computational demands associated with repeatedly solving a new LP can be substantially diminished, if the problem satisfies the non-degenerate condition.

\section{\bfseries\scshape{Non-degenerate Problems and Improved Convergence Rate}}    \label{sec:non-degenerate}

In this section, we define what we call a non-degenerate problem in Section~\ref{subsec:LP-update-policy-and-non-degenerate-property}, and show how it allows one to design an improved LP-update policy with selective updates that is more efficient in Section \ref{subsec:improved-LP-update}. As we will see, when a problem is non-degenerate, the solutions to the LP starting from initial conditions $\bx\approx\bx^*(t)$ can be expressed as linear functions of $\bx$. We will show that this can be used to improve both the computational efficiency of the algorithm and the rate at which the algorithm becomes asymptotically optimal: we prove that the new LP-update policy has an $\calO(1/N)$ performance guarantee for non-degenerate problems. We discuss questions related to rounding in Section \ref{subsec:rounding-problem-discussion}, which is an important step when applying the policy. We provide an alternative formulation of non-degeneracy in Section \ref{subsec:alternative-formulation-of-non-degeneracy}, and compare it with various existing notions in the literature. The proofs of the technical results are postponed to Section~\ref{sec:proof-of-main-theorem}.

\subsection{The non-degenerate property} \label{subsec:LP-update-policy-and-non-degenerate-property}
To ease the notation, we denote by $u := d(A+1)$. We start by remarking that the linear program \eqref{eq:relaxed_problem_general} can be decomposed recursively: By Bellman's principle of optimality, for each $0 \le t \le T-1$, we have
\begin{maxi!}|s|{\by\in\R^u}{\rr^\top \by + \Vrel(\phi(\by),T-t-1) \label{eq:define-argmax} }{\label{eq:one-step-relation-of-V-rel}}{\Vrel(\bx,T-t)=}
  \addConstraint{ \ \ \by \ge \mathbf{0} \label{eq:positive4}}
  \addConstraint{ D \by \le \bb \label{eq:relaxed-constraints-1}}
  \addConstraint{ E \by = \bx \label{eq:init_general-1}}.
\end{maxi!}
where $E$ is the matrix encoding the equality constraints $\sum_{a \in \calA} y_{s,a} = x_s$, so that it has $d$ rows and $u$ columns, for which the $s$-th row is such that from the $(s-1)(A+1)+1$-th to the $s(A+1)$-th coordinates are $1$, the rest are $0$.

Let $\by^*$ be an optimal solution to the linear program \eqref{eq:relaxed_problem_general}, with the corresponding configuration vectors $\bx^*$. Write $\by^*(t)$ (resp. $\bx^*(t)$) the part of $\by$ (resp. $\bx^*$) at time-step $t$. Define $\calJ^*(t)$ as the set of indices for which the budget constraint \eqref{eq:relaxed-constraints-1} is an equality: $(D\by^*(t))_j = b_j$ for all $j\in\calJ^*(t)$ and $(D\by^*(t))_j < b_j$ for $j\in\calJ\setminus\calJ^*(t)$. We now consider the following optimization problem:
\begin{maxi!}|s|{\by\in\R^u}{\rr^\top \by + V_{\mathrm{rel}}(\phi(\by),T-t-1) }{\label{eq:opti_f}}{F_{\by^*}(\bx,T-t)=}
  \addConstraint{ \ \ \by \ge \mathbf{0} \label{eq:relaxed-constraints->0}}
  \addConstraint{ (D \by)_j < b_j\qquad \forall j\not\in\calJ^*(t) \label{eq:relaxed-constraints-<}}
  \addConstraint{ (D \by)_j = b_j\qquad \forall j\in\calJ^*(t) \label{eq:relaxed-constraints-=}}
  \addConstraint{ E \by = \bx \label{eq:relaxed-constraints-=2}}
\end{maxi!}
As \eqref{eq:opti_f} is more restrictive than \eqref{eq:one-step-relation-of-V-rel}, we have $F_{\by^*}(\bx,T-1)\le V_{\mathrm{rel}}(\bx,T-t)$. Moreover, by definition, when $\bx=\bx^*(t)$, we have $F_{\by^*}(\bx^*(t),T-1)=V_{\mathrm{rel}}(\bx^*(t),T-t)$. In what follows, we look for a condition which guarantees that this equality is preserved in a neighbourhood of $\bx^*(t)$, for all $1 \le t \le T-1$. In other words, the sets of saturated and non-saturated constraints of an optimal LP solution to \eqref{eq:relaxed_problem_general} are unchanged upon a small perturbation to the state vectors $\bx^*(t)$.

To do so, let $\calU^*(t) \subset \calU$ be the set of indices $(s,a)$ for which $y^*_{s,a}(t) = 0$, and $\calS^*(t)$ be the set of states for which $x_s^*(t)>0$. We collect all the equality constraints from \eqref{eq:relaxed-constraints->0}, \eqref{eq:relaxed-constraints-=} and \eqref{eq:relaxed-constraints-=2} below:
\begin{align}
  y_{s,a} &= 0 & \forall (s,a)\in\calU^*(t)\label{eq:equalityC0}\\
  \sum_{s,a} D_j(s,a)y_{s,a} &= b_j & \forall j\in\calJ^*(t)   \label{eq:equalityC1}  \\
  \sum_{a} E(s,a)y_{s,a} &= x_s(t) & \forall s\in\calS^*(t)   \label{eq:equalityC2}
\end{align}
The above equalities \eqref{eq:equalityC0}--\eqref{eq:equalityC2} can be represented in compact form using a matrix $C^*(t)$ that has $|\calJ^*(t)| + |\calS^*(t)|+|\calU^*(t)|$ rows and $u=|\calU|$ columns:
\begin{equation} \label{eq:linear-relation-for-y}
  C^*(t) \by = \left[ \mathbf{0}|_{\calU^*(t)}; \ \bb |_{\calJ^*(t)}; \ \bx |_{\calS^*(t)} \right]^\top,
\end{equation}
where the notations $\bb|_{\calJ^*(t)}$ and $\bx|_{\calS^*(t)}$ indicate that the vectors are restricted to the coordinates $\calJ^*(t)$ or $\calS^*(t)$. We wish to establish an inverse to the above linear relation \eqref{eq:linear-relation-for-y}, so that if $\bx^*(t)$ is slightly perturbed as $\widetilde{\bx}$ (where the perturbation comes from the stochastic noise), the inverse equation gives a $\widetilde{\by}$ that remains optimal.

We are prepared to introduce the concept of non-degeneracy, a pivotal term in our discussion.

\begin{defi}[Non-degeneracy] \ \label{def:non-degeneracy}
  An LP \eqref{eq:relaxed_problem_general} is \emph{non-degenerate} if there exists a solution $\by^*$ to the LP, for which at \emph{all} time-step $t\in\{1\dots T-1\}$, the matrix $C^*(t)$ in \eqref{eq:linear-relation-for-y} has rank $|\calJ^*(t)| + |\calS^*(t)|+|\calU^*(t)|$.
\end{defi}
Note that in the above definition we exclude $t=0$, since the perturbation comes from stochastic noise, which only appears from $t=1$ onward.

Recall that our model is a generalization of the restless bandit model studied in \cite{gast2023linear,zhang2021restless,brown2023fluid}, for which a notion of non-degeneracy is already introduced in the aforementioned papers. We show in Section~\ref{subsec:alternative-formulation-of-non-degeneracy} that our definition coincides with their definition when we restrict to the model studied in \cite{gast2023linear,zhang2021restless,brown2023fluid}. This shows that our new definition is indeed an extension of theirs to the broader class of multi-action multi-constraint bandit problems.

\begin{figure}[ht]
  \centering  
  \begin{tabular}{cc}
    \includegraphics[width=0.45\linewidth]{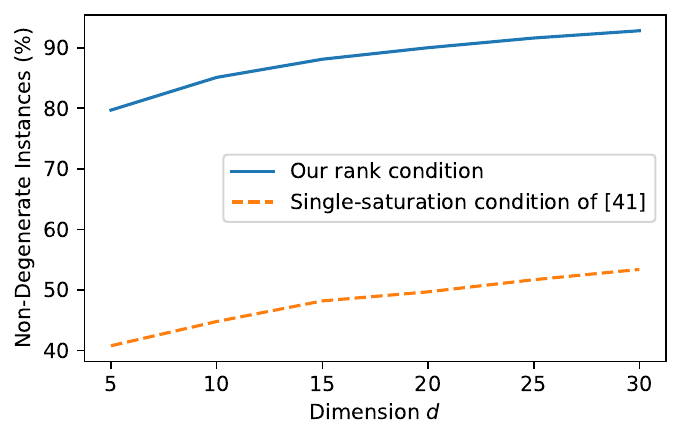}
    &\includegraphics[width=0.45\linewidth]{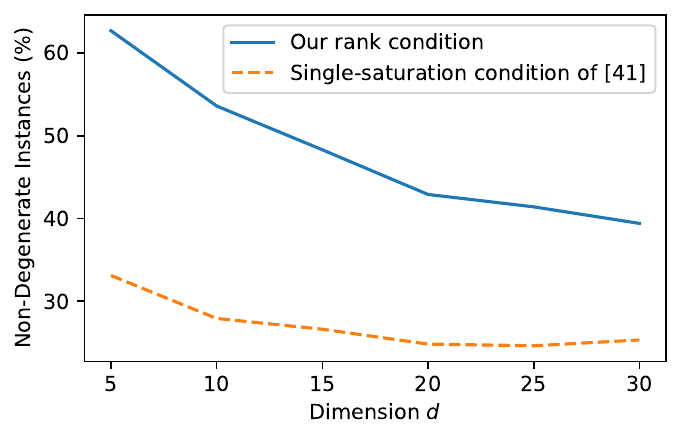}\\
    (a) Transition matrices are dense
    &(b) Transition matrices are tri-diagonal
  \end{tabular}
  \caption{Percentage of non-degenerate examples as a function of the sub-component's state space dimension $d$.}
  \label{fig:comparison_notions_degeneracy}
\end{figure}

To the best of our knowledge, the most comparable notion of non-degeneracy at the same level of generality as ours is found in \citet{zhang2022near1}. However, as discussed in Section~\ref{subsec:alternative-formulation-of-non-degeneracy}, their definition is more restrictive: basically, it asks that a single constraint is saturated. Therefore, our definition encompasses a broader spectrum of problems as non-degenerate compared to theirs. This distinction is visually represented in Figure~\ref{fig:comparison_notions_degeneracy}, where we quantified the proportion of non-degenerate problems according to both definitions. For this, we used randomly generated transition matrices (both dense and tri-diagonal) and reward vectors, plotting the percentage of non-degenerate examples against the number of states $d$. These plots reveal that about half of the problems deemed non-degenerate by our definition are classified as degenerate by the definition in \cite{zhang2022near1}.

Our concept of non-degeneracy, defined through a rank condition, bears a significant resemblance to the Linear Independence Constraint Qualification found in general nonlinear optimization problems, as discussed in works such as \citet{bazaraa2013nonlinear}. The latter stipulates that the gradients of active constraints in an optimization problem must be linearly independent. We adopt this principle in a unique way, by breaking down the linear program outlined in \eqref{eq:relaxed_problem_general} into a $T$-step optimization problem, following Bellman's principle of optimality. In our approach, the notion of non-degeneracy implies that the saturated constraints in the optimal solution of the LP should be distributed "evenly" across the $T$ steps. As we shall see, this property is key to ensuring that we have sufficient rank to create an inverse linear mapping around the original optimal solution at \emph{any} given time-step $t$, while accounting for the stochastic uncertainties in the system states $\bx(t)$.

\subsection{The improved LP-update policy}
\label{subsec:improved-LP-update}

For a given $\bx\in\Delta^d$ and $\varepsilon>0$, we define the neighbourhood of $\bx$ of size $\varepsilon$ as $\calB(\bx,\varepsilon):=\left\{\bx'\in\Delta^d \mid |x'_s-x_s| < \varepsilon \text{ and $x'_s=0$ for all $s$ such that $x_s=0$}\right\}$.  As we show below, being non-degenerate implies that the optimal solution $\by(\bx)$ to the LP \eqref{eq:one-step-relation-of-V-rel} can be chosen as a \emph{locally linear} function of $\bx$ in a small neighbourhood of $\bx^*$.
\begin{prop}
  \label{prop:optimal} \
  Assume that the problem is non-degenerate as in Definition \ref{def:non-degeneracy}. Then, for all time-step $t\ge1$, the matrix $C^*(t)$ has a right inverse $C^+(t)$. Moreover, there exists $\varepsilon>0$ such that:
  \begin{itemize}
    \item The function $\bx\mapsto\relm{\bx,T-t}$ is linear on $\calB(\bx^*(t),\varepsilon)$.
    \item Choosing
      \begin{equation}  \label{eq:linearDef}
      \by(\bx) := \by^*(t) + C^+(t)\left[\begin{array}{c}
      \mathbf{0}|_{\calU^*(t)}\\
      \mathbf{0}|_{\calJ^*(t)}\\
      (\bx-\bx^*(t)) |_{\calS^*(t)}\end{array}\right]
       \end{equation}
    is an optimal solution of \eqref{eq:one-step-relation-of-V-rel} for all $\bx\in\calB(\bx^*(t),\varepsilon)$.
  \end{itemize}
\end{prop}

\begin{proof}{Proof.}
  By standard linear algebra argument, a matrix of dimension $d_1\times d_2$ with rank $d_1$ has a right inverse. This implies that there exists a matrix $C^+(t)$ such that $C^*(t)C^+(t)$ is the identity matrix. Hence if $\by(\bx)$ is defined by \eqref{eq:linearDef}, then
  \begin{align*}
    C^*(t)\by(\bx) = C^*(t)\by^*(t) + C^*(t)C^+(t)\left[\begin{array}{c}
      \mathbf{0}|_{\calU^*(t)}\\
      \mathbf{0}|_{\calJ^*(t)}\\
      (\bx-\bx^*(t)) |_{\calS^*(t)}
    \end{array}\right]
    = \left[\begin{array}{c}
      \mathbf{0}|_{\calU^*(t)}\\
      \mathbf{b}|_{\calJ^*(t)}\\
      \bx |_{\calS^*(t)}
    \end{array}\right]
  \end{align*}
  In particular, $\by(\bx)$ satisfies \eqref{eq:equalityC0}--\eqref{eq:equalityC2}. This shows that there exists $\varepsilon>0$ such that $\by(\bx)\in\calY(\bx)$ (i.e. satisfy all constraints of \eqref{eq:one-step-relation-of-V-rel}) for all $\bx\in\calB(\bx^*(t),\varepsilon)$, because the constraints that are not covered by \eqref{eq:equalityC0}--\eqref{eq:equalityC2} are either satisfied by \emph{strict} inequalities for $\bx^*(t)$, or correspond to $x^*_s(t)=0$. It is worth noting that the selection of $\varepsilon$ is explicitly based on the \emph{slack} of the non-saturated constraints.

  We now prove by a backward induction on $t$, that for all $t$, there exists $\varepsilon_t>0$ such that the function $\bx\mapsto\relm{\bx,T-t}$ is linear on $\calB(\bx^*(t), \varepsilon_t)$.  This is clearly true for $t=T$ for which $\relm{\bx,0}=0$ for all $\bx$.

  Assume now that it holds for some $t+1\le T$, and denote by $g_t(\bx,\by)$ the reward provided by the control $\by$ with initial state vector $\bx$ and horizon $T-t$. As shown before, for $\bx$ close enough to $\bx^*(t)$, the control $\by(\bx)$ is feasible for $\bx$. Moreover, the induction hypothesis implies that $\relm{\phi(\by(\bx)),T-t-1}$ is locally linear in $\bx$ for all $\bx$ close enough to $\bx^*(t)$. This shows that $\bx\mapsto g_t(\bx,\by(\bx))$  is locally linear on $\calB(\bx^*(t),\varepsilon_t)$.  We argue that this implies that $\by(\bx)$ is the optimal control for all $\bx\in\calB(\bx^*(t),\varepsilon_t)$. Indeed:
  \begin{itemize}
    \item As $V_{\mathrm{rel}}(\bx,T-t)$ is the solution of a linear program where $\bx$ can be seen as a linear constraint, the function $\bx \mapsto V_{\mathrm{rel}}(\bx,T-t)$ is concave in $\bx$. Moreover, by construction $\by(\bx^*(t))$ provides the optimal solution for $\bx^*(t)$.
    \item Denote by $\bx' := 2\bx^*(t) - \bx$. By construction $\bx'\in\calB(\bx^*(t),\varepsilon_t)$. By concavity of the $V_{\mathrm{rel}}(\cdot, T-t)$ function on $\bx$, the possible sub-optimality of $\by(\bx)$, and together with the local linearity of $g_t(\bx,\by(\bx))$, we deduce that:
    \begin{align*}
      V_{\mathrm{rel}}(\bx^*,T-t) &\ge (V_{\mathrm{rel}}(\bx,T-t)+V_{\mathrm{rel}}(\bx',T-t))/2 \\
      &\ge (g_t(\bx',\by(\bx')) + g_t(\bx,\by(\bx))/2 \\
      & = g_t(\bx^*,\by(\bx^*)) = V_{\mathrm{rel}}(\bx^*,T-t).
    \end{align*}
    This shows that the inequalities must be equalities, which implies that $\by(\bx)$ defined as in Equation \eqref{eq:linearDef} is optimal on $\calB(\bx^*(t),\varepsilon_t)$.    \qed
  \end{itemize}

\end{proof}


The definition of a non-degenerate problem suggests that the original LP-update Algorithm~\ref{algo:LP-update} can be implemented by only recomputing an update (i.e. re-solving) when necessary:
\begin{enumerate}
  \item \label{update:case 1} When at a time $t$ the rank condition on $C^*(t)$ is not satisfied, then one cannot compute the right inverse $C^+(t)$.
  \item \label{update:case 2} When at a time $t$, the rank condition is satisfied but the stochastic trajectory has deviated too much from the optimal deterministic trajectory, so that the suggested decision vector $\by(\bXN(t))$ in \eqref{eq:linearDef} no longer gives a feasible decision vector, i.e. is not in $\calY(\bXN(t))$.
\end{enumerate}
When one of these two situations occurs, the algorithm solves a new LP with initial configuration vector $\bXN(t)$ over the horizon $[t,T]$.  This leads to the improved LP update algorithm that is described in Algorithm~\ref{algo:LP-update-improved}.


\begin{algorithm}[t]
\SetAlgoLined
\SetKwInput{KwInput}{Input}
\KwInput{Initial configuration vector $\bXN(0)$.}
  \For{$t=0,1,2,\dots,T-1$}{
    \eIf{$t\ge1$ and $C^*(t)$ has rank $|\calJ^*(t)| + |\calS^*(t)|+|\calU^*(t)|$ and $\by(\bXN(t))\in\calY(\bXN(t))$}{
      Set $\by(t) := \by(\bXN(t))$ as in \eqref{eq:linearDef}\;
    }
    {
      Solve LP \eqref{eq:relaxed_problem_general} with parameters $(\bXN(t), T-t)$ and call the output $\by^*(\tau)$ for $\tau\in\{t\dots T\}$\;
      Set  $\by(t) := \by^*(t)$\;
    }
    Set  $\YN_{s,a}(t):=N^{-1}\floor{Ny_{s,a}(t)}$ for $a\not= 0$ and  $\YN_{s,0}(t) := \XN_{s}(t) - \sum_{a\not = 0} \YN_{s,a}(t)$ \;
  Use actions $\YN_{s,a}(t)$ over all sub-MDPs to advance to the next time-step\;
 }
 \caption{The LP-update policy for weakly-coupled MDPs (selective updates).}
 \label{algo:LP-update-improved}
\end{algorithm}

It is worth noting that $\by(t)\in\calY(\bXN(t))$ is not a sufficient condition for optimality (this only holds if $\bXN(t)$ is in the neighborhood of $x^*(t)$).  Consequently, Algorithm~\ref{algo:LP-update-improved} may not have the same output as Algorithm~\ref{algo:LP-update}, even in cases where no update is performed in Algorithm~\ref{algo:LP-update-improved}. Nonetheless, the square-root convergence rate established in Theorem \ref{thm:LP-update-general} is also applicable to Algorithm~\ref{algo:LP-update-improved}, as can be demonstrated using analogous reasoning.

Additionally, it is important to highlight that Algorithm~\ref{algo:LP-update-improved} requires fewer updates than Algorithm~\ref{algo:LP-update}, enhancing its computational efficiency, see Section \ref{subsec:computation-time-analysis} for a demonstration. In our forthcoming proof, we will show that the discrepancy in the outcomes of the two algorithms is exponentially small for non-degenerate problems. Therefore, in the subsequent theorem, we use $\Update$ to refer interchangeably to the value derived from either the LP-update policy in Algorithm~\ref{algo:LP-update} or Algorithm~\ref{algo:LP-update-improved}.

\subsection{Performance guarantees}

The next theorem shows that when a problem is non-degenerate, the sub-optimality of the LP-update is of order $\calO(1/N)$, which is an order of magnitude faster than the $\calO(1/\sqrt{N})$ guaranty of Theorem~\ref{thm:LP-update-general}. One drawback, though, is that the hidden dependence on $T$ in the $\calO(1/N)$ can be exponential in $T$. This exponential dependence on $T$ comes from the proof of the theorem which can be decomposed in two steps: first, we show that $\bXN(t)$ is close to $\bx^*(t)$ that when $N$ is large, and second we show that $\Vrel$ is local linear close to $\bx^*(t)$. The first step leads to an exponential dependence in $T$ in the worst case. This dependence disappears in the proof of Theorem~\ref{thm:LP-update-rate-nondegenerate} because we only use that $\Vrel$ is Lipschitz-continuous but not that $\bXN(t)$ is close to $\bx^*(t)$. That being said, in practice, the LP-update becomes asymptotically much faster for non-degenerate examples than for degenerate examples (see Section~\ref{ssec:experiment-nondegeVSdege}).

\begin{thm} \label{thm:LP-update-rate-nondegenerate} \
  Denote by $\Update$ the value of the LP-update policy defined in Algorithm~ \ref{algo:LP-update-improved} (or  Algorithm~\ref{algo:LP-update}), and by $\rel$ the value of the linear program \eqref{eq:relaxed_problem_general}.
  \begin{itemize}
    \item[(i)] For any non-degenerate weakly-coupled MDP with statistically identical arms, there exists constant $C>0$ such that for all $N$:
    \begin{equation*}
      \abs{\Update - \rel} \le \frac{C}{N}.
    \end{equation*}
    \item[(ii)] There exists a non-degenerate weakly-coupled MDP with statistically identical arms and a constant $C'>0$ such that for all $N$:
    \begin{equation*}
      \abs{\Update - \rel} \ge \frac{C'}{N}.
    \end{equation*}
 \end{itemize}
\end{thm}
\begin{proof}{Proof.}
  As for Theorem~\ref{thm:LP-update-general}, we prove the upper bound (i) and postpone the proof of the lower bound (ii) to Section~\ref{subsec:lowerbound}.

  Let $\by^*$ be the optimal solution of the LP solved at time $0$ and let $x^*_s(t) = \sum_a y^{*}_{s,a}(t)$. In what follows, we show that with very high probability, $\bXN(t)$ remains close to $\bx^*(t)$ so that no update is necessary, and the sequence of $\bY(t)$ computed by Algorithm~\ref{algo:LP-update-improved} will always be an optimal solution to the LP problems. To this end, we shall bound the probability of the event when there is some $t$ for which $\bXN(t)$ deviates far away from $\bx^*(t)$ from above.

  By Proposition~\ref{prop:optimal}, there exists $\varepsilon>0$ such that $\relm{\bx,T-t}$ is linear in $\calB(\bx^*(t),\varepsilon)$ and such that the control $\by(\bx)$ defined in Proposition~\ref{prop:optimal} and used in Algorithm~\ref{algo:LP-update-improved} is optimal when $\bXN(t) \in \calB(\bx^*(t),\varepsilon)$. Call $\calE$ the event $\bXN(t)\in\calB(\bx^*(t),\varepsilon)$. By Lemma~\ref{lem:concentration_trajectory}, there exists $C_1,C_2>0$ such that the event $\calE$ holds with probability at least $C_1 \cdot e^{-C_2N}$.  Hence, when $\calE$ is true, Algorithm~\ref{algo:LP-update-improved} behaves as Algorithm~\ref{algo:LP-update}. This shows that \eqref{eq:proof_update} holds also for Algorithm~\ref{algo:LP-update-improved}, up to an $\calO(e^{-C_2N})$ term due to the (exponentially small) probability that $\calE$ does not hold. Similarly \eqref{eq:Z(t)} also holds for Algorithm~\ref{algo:LP-update-improved} up to an $\calO(e^{-C_2N})$ term. As $\relm{\bx,T-t}$ is locally linear when $\calE$ holds, \eqref{term:B} of \eqref{eq:Z(t)} is smaller than $C_1 \cdot e^{-C_2N}$.

  This shows that the sub-optimality gap is $\sum_{t=0}^{T} \rr^\top \expect{\bYN(t)-\bY(t)} + \calO(e^{-C_2N})$. Moreover, by our explicit construction of $\bYN(t)$ as in \eqref{eq:rounding}, it satisfies $\norme{\bYN(t)-\bY(t)} \le d|\calA|/N$, which gives Theorem~\ref{thm:LP-update-rate-nondegenerate}.
  \qed
\end{proof}

It is noteworthy that the term $\calO(1/N)$ in our analysis arises totally from the rounding procedure detailed in \eqref{eq:rounding}. This aspect is further explored in Section~\ref{subsec:rounding-problem-discussion}, where we introduce the concept of a \emph{perfect rounding condition}. Under this ideal condition, the rounding error is effectively eliminated, thereby enhancing the convergence rate to $\calO(e^{-C_2N})$. Additionally, this section presents an alternative rounding strategy, which has demonstrated superior performance in practical applications.


\subsection{Sub-optimality gap for degenerate and non-degenerate problems}  \label{ssec:experiment-nondegeVSdege}

\begin{figure}[ht]
  \centering
  \begin{tabular}{cc}
    \includegraphics[width=0.5\linewidth]{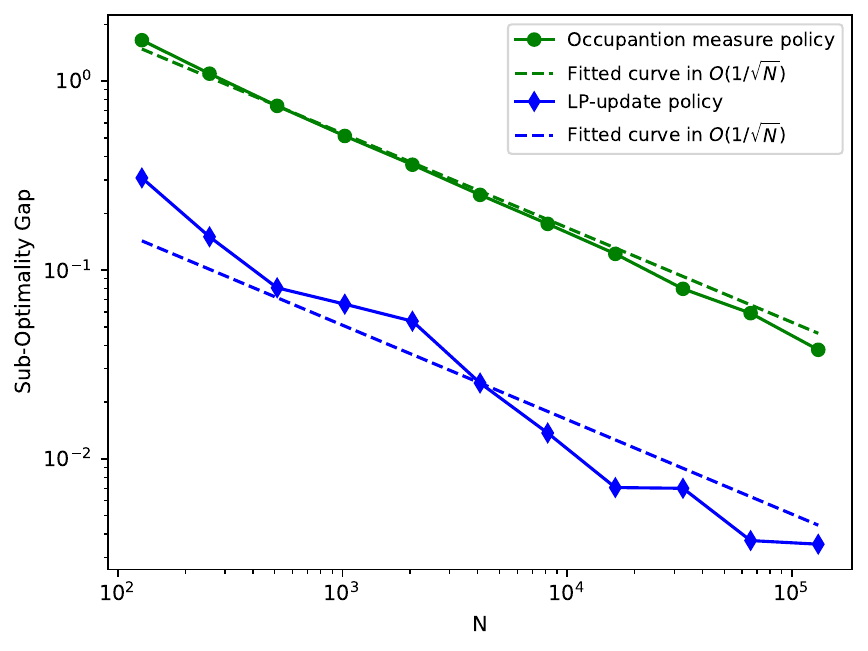}
    &\includegraphics[width=0.5\linewidth]{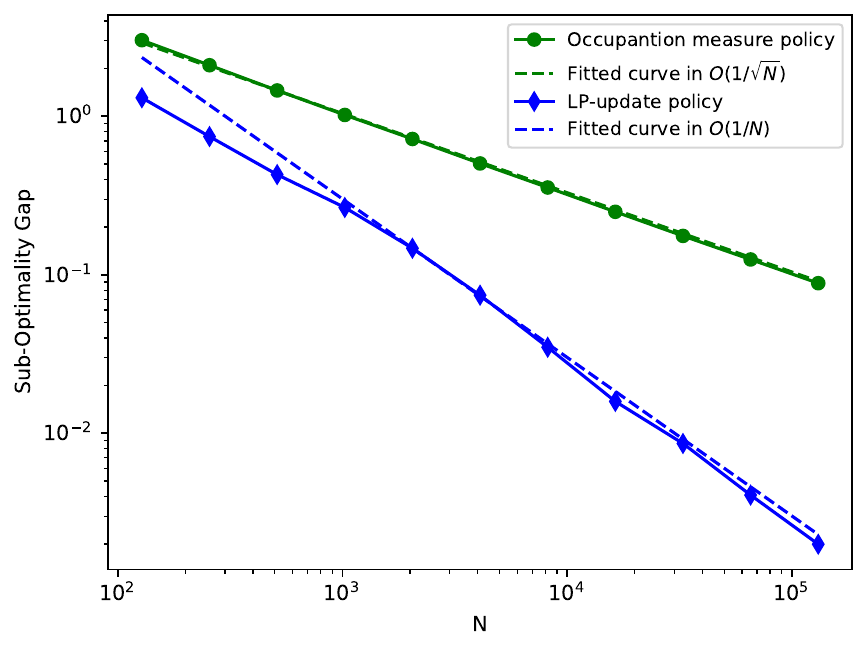}
    \\
    (a) A degenerate example
    &(b) A non-degenerate example
  \end{tabular}
  \caption{Sub-optimality gap in degenerate and non-degenerate examples: for non-degenerate examples, the sub-optimality gap of our LP-update policy decreases at a faster rate $\calO(1/N)$ than the occupation measure policy of \cite{ZayasCabn2017AnAO}.
  }
  \label{fig:dege_vs_nondege}
\end{figure}

To demonstrate the distinct behaviors between degenerate and non-degenerate cases, we conducted an experiment with two randomly chosen examples – one degenerate and one non-degenerate. The results, showcased in Figure~\ref{fig:dege_vs_nondege}, compare the sub-optimality gap of the LP-update policy, as well as the occupation measure policy as a function of the number of arms, $N$. This comparison unequivocally shows that our LP-update policy consistently outperforms the others.

Additionally, the experiment highlights notable differences in the rate at which these policies approach asymptotic optimality. Consistent with our theorem, the sub-optimality gap for the LP-update policy is $\calO(1/\sqrt{N})$ in degenerate cases and improves to $\calO(1/N)$ in non-degenerate cases. In contrast, the occupation measure policy suffers a sub-optimality gap of $\calO(1/\sqrt{N})$ in the two scenarios. 

We remark that the zigzag pattern in the performance curves of the LP-update policy is caused by cutoff rounding errors, as outlined in Equation~\eqref{eq:rounding}. This pattern is mitigated when a more sophisticated rounding approach is employed, as showcased in Section \ref{subsec:sophisticated-rounding}.

\section{\bfseries\scshape{Case Study: Generalized Applicant Screening Problem}} \label{sec:general-applicant-screening}

The applicant screening problem is proposed in \citet{Brown2020IndexPA}, and has been subsequently studied in \citet{gast2023linear}. The problem studied in these papers can be modeled as a two-action single-constraint restless bandit problem. We study in this section a generalization of the problem allowing multiple actions while also having multiple constraints, so that the model fits naturally into our weakly-coupled MDP framework.


\subsection{Problem description} \label{subsec:problem-description}

Consider a group of $N$ applicants applying for a job. The decision maker's goal is to hire the best possible $\beta N$ applicants. The applicant $n$ has an unknown quality level $p_n \in [0,1]$. At each decision epoch $t$, the decision maker chooses for each candidate either to interview this candidate with one or two questions, or chooses not to interview this candidate, giving the action set $A = \{0,1,2 \}$. For each interviewed applicant (i.e. $a \in \{1,2\}$), a signal $q_n(t) \sim \text{binomial}(a,p_n)$ is returned, indicating how many among the $a$ questions have been solved correctly by the applicant. All variables $q_n(t)$ are supposed to be independent.

Choosing an action $a$ on an applicant consumes $D(a)$ units of resource (time, space, organization cost, etc.), for which we choose to be
\begin{equation} \label{eq:resource-consumption}
  D(s,a) = D(a) = \begin{cases}
  0, & \mbox{if } a=0; \\
  1, & \mbox{if } a=1; \\
  1.5, & \mbox{if } a=2.
  \end{cases}
\end{equation}
The value "$1.5$" on action $a=2$ is interpreted as asking a single applicant consecutively two questions consumes less resource than asking separately two applicants each with one question. At each decision epoch a total amount of $\alpha N$ resource is available. There is a total number of $T$ interviewing rounds, and in the final $(T+1)$-th round, the decision maker admits $\beta N$ applicants, based on the results of the interviewing rounds.

We assume that the applicants belong to two different groups, each having a population of $N/2$. And we will consider two scenarios:
\begin{itemize}
  \item No fairness: in this scenario, there is a single budget constraint of $\alpha N$ for the whole population.
  \item Fair selection: in addition to the above constraint, the decision maker cannot spend more that $\gamma N$ budget on each of the single group, where $\gamma < \alpha < 2 \gamma$.
\end{itemize}

Notice that the above applicant screening problem generalizes the problem studied in \citet{Brown2020IndexPA} and \citet{gast2023linear} in two ways: we allow for more than one question whereas the cited paper consider $A=\{0,1\}$; and we consider a fairness constraint, whereas the cited paper only have a single resource constraint that the number of interviewed candidate should not be larger than $\alpha N$.

\subsection{Modeling as weakly-coupled MDP} \label{subsec:modeling-the-problem}

To cast the problem into a weakly-coupled MDP as described in Section \ref{sec:model-description}, we consider a Bayesian model in which the quality level $p$ of an applicant from each group is generated from some beta distribution, and the decision maker's estimation on each applicant's $p$ is updated using Bayes' rule. The state $s$ of an applicant is hence a $2$-tuple $(a,b)$, indicating the current estimation of her quality level and is distributed according to the beta distribution $beta(a,b)$.

For the first $T$ interviewing rounds, the action set is $A = \{0,1,2 \}$. Upon taking action $0$ on an applicant, the estimation is unchanged, so the matrix $\pp^0$ is the identity matrix. Upon taking action $1$, the state is updated according to
\begin{equation*}
  (a,b) \xrightarrow[]{\text{action $1$}} \begin{cases}
                                              (a+1,b), & \mbox{with probability $a/(a+b)$; }  \\
                                              (a,b+1), & \mbox{with probability $b/(a+b)$. }
                                            \end{cases}
\end{equation*}
This gives the matrix $\pp^1$. Likewise, for action $2$ the state is updated as
\begin{equation*}
  (a,b) \xrightarrow[]{\text{action $2$}} \begin{cases}
                                              (a+2,b), & \mbox{with probability $\frac{a(1+a)}{(a+b)(1+a+b)}$; }  \\
                                              (a+1,b+1), & \mbox{with probability $\frac{2ab}{(a+b)(1+a+b)}$; }  \\
                                              (a,b+2), & \mbox{with probability $\frac{b(1+b)}{(a+b)(1+a+b)}$. }
                                            \end{cases}
\end{equation*}
This gives the matrix $\pp^2$. The function $D(s,a)$ is given by \eqref{eq:resource-consumption} and is independent of the state $s$. The rewards are all zero during the first $T$ rounds.

For the final $(T+1)$-th admitting round, the action set is $\{ \text{admit, not admit} \}$. The reward on an admitted applicant in state $(a,b)$ is $a/(a+b)$; the reward on a non admitted applicant is zero.

\subsection{Discussion on simulation results} \label{subsec:simulation-results}

For our simulations, we choose $\beta = 0.1$, and we consider two scenarios. The first one is such that $\alpha = 0.15$, $\gamma = 0.1$. This is the scenario where the resource is "scarce". The second scenario where the resource is "abundant" is such that $\alpha=0.3$ and $\gamma=0.2$ (resource being doubled). In each scenario we shall compare the effect of adding or removing the fairness constraints. Without the fairness constraints, the decision maker can distribute the total $\alpha N$ units of resource freely  among the two groups at  each interviewing round.

We assume that the decision maker's prior estimations on the two groups of applicants are respectively $beta(1,1)$ and $beta(2,2)$, so that they have the same mean but the second group has a lower variance. Throughout the horizon is fixed to  $T=10$. To make the problem more realistic, and in the same time to reduce the total number of possible states of the MDP, we require in addition that no more than $10$ questions can be asked on a single applicant during the $10$ interviewing rounds.

\begin{figure}[hbtp]
  \centering
    \includegraphics[width=1.\linewidth]{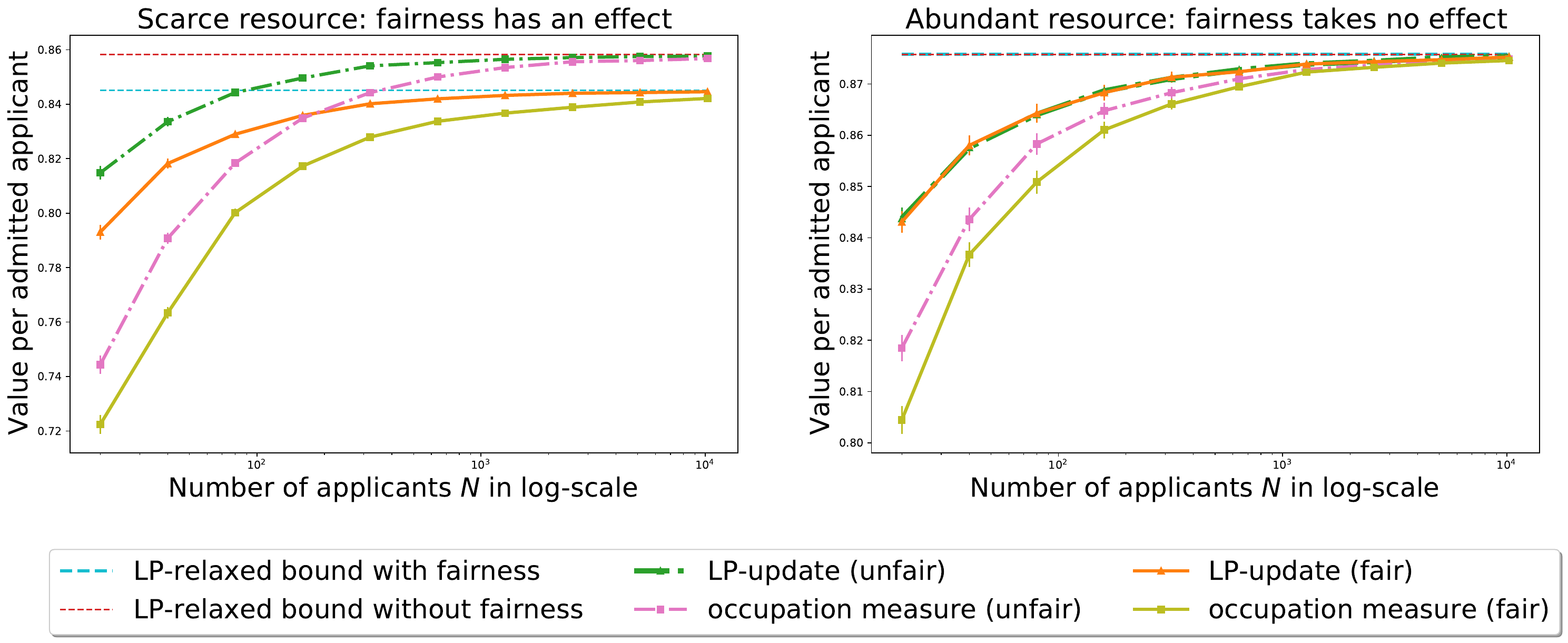}
    \caption{Performance on generalized applicant screening problems when the resource is scarce (left panel) or abundant (right panel), with or without fairness constraints. Overall the LP-update policy outperforms the occupation measure policy in all situations, and the advantage is more apparent for small values of $N$. Fairness has a negative impact in the scarce resource scenario, whereas in the abundant case, it has no effect asymptotically, but can still influence the performance of the occupation measure policy. }
  \label{fig:performance_generalaized_applicant_screening}
\end{figure}

The simulations are done for $N$ ranging from $20$ to $10240=20\times 2^{10}$. For each value of $N$, we generate $1600$ instances of $N$ applicants according to the beta distributions we described previously in each scenario,. We apply the LP-update policy and occupation measure policy on each instance, with or without the fairness constraints. In each simulation the performance is evaluated based on the average quality level of the final admitted $\beta N$ applicants. The results are reported in Figure \ref{fig:performance_generalaized_applicant_screening}.

Note that as guaranteed by Theorem \ref{thm:LP-update-general} and Theorem $1$ of \citet{xiong2021reinforcement}, both LP-update policy and occupation measure policy converge to the  LP-relaxed bounds, as $N$ goes to infinity. This is what we observe in Figure~\ref{fig:performance_generalaized_applicant_screening}: all policies converge to the relaxed upper bound as $N$ goes to infinity. The situation is however different for smaller values of $N$: In all cases, the LP-update policy outperforms the occupation measure policy. The smaller the value of $N$, the more apparent is the advantage of the LP-update policy.

We shall now discuss the effect of adding fairness constraints. We observe that in the first scarce resource scenario,  the LP-relaxed performance bound becomes smaller when  fairness is imposed. This should be expected since adding more constraints on a linear program can only decrease its optimal value. However,  in the second scenario, when the resource is doubled, adding or removing the fairness constraints result in the same upper-bound value, and the performances of LP-update policy are identical in the two situations. Nevertheless, these fairness constraints still play a role in both policies: for the LP-update policy we observe that the fairness constraints can be saturated when solving new linear programs by applying updates, but it turns out that  this does not influence its performance at all. For the occupation measure policy, however, sampling an action is more restrictive if there are more constraints, and as a result this influences its performance, even though asymptotically they converge to the same limit.

\subsection{Computation time analysis} \label{subsec:computation-time-analysis}

Although the LP-update policy outperforms the occupation measure policy under all situations, it comes with a price of using more computations. This extra computation time is due to the fact that the LP-update policy periodically solves a new LP problem.  We record the execution time for different values of $N$ in Table~\ref{tab:computation-time}, averaging over $100$ runs on one problem instance under each circumstance, together with the $95 \%$ confidence interval. Here by execution time we mean the time needed to apply a policy on an instance of $N$ applicants problem until the final phase of admission round. The program is written in $\mathrm{Python}$ using $\mathrm{NumPy}$ for data structure and $\mathrm{PuLP}$ for solving LP's. Note that our code is not particularly optimized and is tested on an ordinary personal laptop.

\begin{table}[ht]
  \centering
  \begin{tabular}{|ccc|ccc|}
  \hline
    &&\multicolumn{1}{r|} {Number of applicants} & $N=20$ & $N=100$ & $N=1000$  \\ \hline
    \parbox[t]{2mm}{\multirow{3}{*}{\rotatebox[origin=c]{90}{\small With}}}
    &\parbox[t]{2mm}{\multirow{3}{*}{\rotatebox[origin=c]{90}{fairness}}}
    &Occupation measure    & 2.70 $\pm$ 0.01  & 2.72 $\pm$ 0.01   & 2.83 $\pm$ 0.01    \\
    &&LP-update    & 12.69 $\pm$ 0.39  & 10.44 $\pm$ 0.41  & 7.06 $\pm$ 0.19   \\
    &&(Number of updates)   &  ($\approx$ 6.4 ) &  ($\approx$ 5.2)  &  ($\approx$ 3.9 ) \\
    \hline
    \parbox[t]{2mm}{\multirow{3}{*}{\rotatebox[origin=c]{90}{\small Without}}}
    &\parbox[t]{2mm}{\multirow{3}{*}{\rotatebox[origin=c]{90}{fairness}}}
    &Occupation measure    & 2.62 $\pm$ 0.01  & 2.64 $\pm$ 0.01 & 2.75 $\pm$ 0.01   \\
    &&LP-update & 8.95 $\pm$ 0.31 & 7.93 $\pm$ 0.23 & 6.75 $\pm$ 0.16   \\
    &&(Number of updates)   &  ($\approx$ 4.5)  &  ($\approx$ 3.6)  &  ($\approx$ 2.8)  \\\hline
  \end{tabular}
  \caption{\label{tab:computation-time}%
    Computation time (in seconds) of the LP-update policy and the occupation measure policy for different values of $N$, as well as the number of times that the LP-update policy solves a new linear programm (the first initial one not included).}
\end{table}

We notice that the variance of the computation time is much larger for the LP-update policy, since the number of time-steps when it  solves a new LP for an update varies on each run. Overall this variance is decreasing as $N$ becomes large. We also remark that the  LP-update policy takes less time for $N=1000$ compared to $N=20$. This is because for larger $N$ the stochastic trajectory is closer to the deterministic one, consequently it is less likely to perform updates caused by a non-feasible action obtained from \eqref{eq:linearDef} (situation $2$ described before Algorithm \ref{algo:LP-update-improved}). Indeed, we observe from Table \ref{tab:computation-time} that the number of updates is decreasing with $N$.

\section{Conclusion, Discussions and Future Works} \label{sec:conclusion-and-future-work}

In this paper, we study finite horizon weakly-coupled Markov decision processes with symmetric components. We propose a re-solving heuristic that is asymptotically optimal under a population scaling of the model. Furthermore, we identify a non-degenerate condition on the model that can improve both the convergence rate and efficiency of the algorithm.

We conclude the main part of the paper by making the following remarks and pointing out several future research directions: 
\begin{enumerate}
  \item In our numerical experiments, we observe that applying re-solving in the LP-update policy improves the performance, as compared to the one-pass policies that do not. However, under other similar problem settings, it is not always better to re-solve. For instance, in the revenue management problem (e.g. \cite{cooper2002asymptotic, jasin2012re, bumpensanti2020re}), there exists a counter-example in \citet{cooper2002asymptotic} where using re-solving gives strictly worse performance in expectation than not using it. A theoretical analysis on the effect of re-solving is provided in \citet{secomandi2008analysis} for this model, giving sufficient conditions so that re-solving is never detrimental. It is unclear whether some similar sufficient conditions exist in our model. Nevertheless, if the problem is non-degenerate, then as can be seen from the proof of Theorems~\ref{thm:LP-update-rate-nondegenerate} and \ref{thm:LP-update-rate-expo}, it is indifferent whether we re-solve or not in the large $N$ limit. Yet, a complete understanding of the re-solving effect in our case is still open. 
      
  \item There exists various possibilities to extend the weakly-coupled MDPs model considered in this paper. These include:
  \begin{itemize}
    \item We could have additional resource constraints that span multiple time-steps, instead of covering only a single time-step. 
    \item The population of the sub-MDPs may not be fixed by $N$, instead there are systems in which sub-component arrive and depart from the system. A typical example is the case of queueing networks where the sub-components are tasks \cite{Ve2016.6}.
    \item The re-solving nature of the LP-update policy requires that the problem be formulated in finite horizon, or to discounted infinite horizon by truncating the time as in \cite{brown2023fluid,zhang2022near}. One may wish to consider the same general weakly-coupled MDPs under infinite horizon. To the best of our knowledge, an asymptotically optimal policy under this generality is still an open question. 
    \item The reward may not be additive over the sub-components (but still additive across time), since it is common to have a utility function that is concave and non-linear. Linearity plays a crucial role in our methodology, especially in the non-degenerate condition and subsequently the faster convergence rate. It is worthwhile to investigate how far we can go by relaxing this linear condition. 
  \end{itemize}

\end{enumerate}

\bibliographystyle{informs2014}
\bibliography{reference}

\begin{APPENDICES}




\section{\bfseries\scshape{Proof of the Technical Results}} \label{sec:proof-of-main-theorem}

This section is dedicated to the proofs of some technical lemmas that are used to prove Theorem~\ref{thm:LP-update-general}, Theorem~\ref{thm:LP-update-rate-nondegenerate} and Theorem~\ref{thm:LP-update-rate-expo}. We also provide counter-examples for the lower bounds of the three theorems, as well as an example to illustrate the equivalence of the two notions of non-degeneracy. 

\subsection{One-step transition and concentration arguments} \label{subsec:preliminary-result}

Recall that the linear function $\phi$ maps a decision vector $\by$ to a configuration vector $\phi(\by) = \big(\phi_1(\by), \dots, \phi_d(\by) \big) \in \Delta^d$ whose $s$-th component is
\begin{equation} \label{eq:Markovian-transition}
  \phi_s(\by) = \sum_{s', a} y_{s',a} P^a_{s',s}.
\end{equation}
This is the deterministic behavior of the Markov transition at time-step $t$. We claim that:
\begin{lem}  \label{lem:Markovian-transition-analysis} \
  Define the random vector $\bEN(t) := \bXN(t+1) - \phi(\bYN(t))$. We have
  \begin{align}
    \expect{\bEN(t) \mid \bYN(t)} &= \mathbf{0},\label{eq:unbiased}\\
    \expect{\norme{\bEN(t)} \mid \bYN(t)} &\le \frac{\sqrt{d}}{\sqrt{N}},\label{eq:CLT}\\
    \proba{\norme{\bEN(t)} \ge \epsilon \mid \bYN(t)} &\le 2d \cdot e^{-2N \epsilon^2/d^2}.\label{eq:hoeffding}
  \end{align}
\end{lem}

\begin{proof}{Proof of Lemma \ref{lem:Markovian-transition-analysis} }
  For simplicity of notation, let us denote by $\by := \bYN(t)$. There are $Ny_{s,a}$ arms in state $s$ and whose action is $a$, and each of these arms makes a transition to  state $s'$ with probability $P^a_{s,s'}$. This shows that $\XN(t+1)$ can be written as a sum of independent random variables as follows:
\begin{align*}
  \XN_{s'}(t+1) = \frac1N\sum_{s,a} \sum_{i=1}^{Ny_{s,a}} \mathbf{1}_{ \{ \sum_{s''=1}^{s'-1} P^a_{s,s''} \le U_{s,a,i} < \sum_{s''=1}^{s'} P^a_{s,s''} \}},
\end{align*}
where $U_{s,a,i}$ are i.i.d uniform random variables in $[0,1]$. Taking expectation then gives $\expect{\XN_{s'}(t+1) \mid \bYN(t)}=\phi_{s'}(\bYN(t))$, which gives \eqref{eq:unbiased}. It also implies that
\begin{align*}
  \expect{|E^{(N)}_{s'}(t+1)|^2 \mid \bYN(t)=\by} &= \var{\XN_{s'}(t+1)\mid \bYN(t)=\by}\\
  & = \frac1{N^2}\sum_{s,a}Ny_{s,a} P^a_{s,s'}(1-P^a_{s,s'})\\
  &\le \frac{\sum_{s,a}y_{s,a} P^a_{s,s'}}{N}.
\end{align*}
This shows that
\begin{align*}
  \expect{\norme{\bEN(t+1)} \mid \bYN(t)=\by} \le  \sqrt{d} \frac{ \sqrt{\sum_{s'} \sum_{s,a}y_{s,a}P^a_{s,s'}}}{\sqrt{N}} = \frac{\sqrt{d}}{\sqrt{N}},
\end{align*}
where the first inequality comes from Cauchy-Schwartz, and this gives \eqref{eq:CLT}.

Equation~\eqref{eq:hoeffding} is a direct consequence of Hoeffding's inequality. Indeed, one has
\begin{equation*}
  \proba{|E^{(N)}_s(t)|\ge\varepsilon/d \mid \bYN(t)} \le 2e^{-N\varepsilon^2/d^2}.
\end{equation*}
By using the union bound, this implies that
\begin{equation*}
  \proba{\norme{\bEN(t)} \ge \varepsilon \mid \bYN(t)}\le d \cdot \proba{|E^{(N)}_s(t)|\ge\varepsilon/d \mid \bYN(t) } \le 2d \cdot e^{-N\varepsilon^2/d^2}.
\end{equation*}
\qed
\end{proof}

\subsection{Lipschitz-constant of $\Vrel(\cdot, t)$}   \label{subsec:lipschitz-const}

Recall that $\gamma$ is equal to 
\begin{align}
  \label{eq:gamma}
    \gamma = \min_{i \neq j, a} \sum_k \min( P^0_{ik}, P^a_{jk}).
\end{align}

\begin{lem}
    \label{lem:Lipschitz-lemma}
    Recall that $\Rmax=\max_{s,a}\abs{r^a_s}$. Then, for all $\bx,\tilde{\bx}\in\Delta$, the value function satisfies: 
    \begin{align*}
        |\Vrel(\bx,t)-\Vrel(\tilde{\bx},t)| \le (1+(1-\gamma) + \dots + (1-\gamma)^{t-1})\Rmax \norme{\bx-\tilde{\bx}}.
    \end{align*}
\end{lem}

In particular, this result implies that:
\begin{itemize}
  \item The Lipschitz constant of $\bx\mapsto\Vrel(\bx,t)$ is upper bounded by $t\Rmax$;
  \item When $\gamma>0$, the Lipschitz constant of $\bx\mapsto\Vrel(\bx,t)$ is bounded by $\Rmax /\gamma$ (for all $t$).
\end{itemize}

\begin{proof}{Proof of Lemma \ref{lem:Lipschitz-lemma} }
     The main idea is to work with vectors of $N$ components $\bs\in\calS$ instead of the population state $\bx$. The use of vectors of components will allow the usage of a sort of coupling between trajectories. To simplify our notation, we will need a way to go from vectors of $N$ components to population vectors. For a population vector $\bs\in\calS^N$, we denote by $\xs(\bs)$ the vector such that for $i\in\calS$: 
    \begin{align*}
        \xs_i(\bs) = \frac1N \sum_{n=1}^N \mathbf{1}_{\{s_n=i\}}.
    \end{align*}
    Recall that $\Delta$ is the simplex of dimension $d=|\calS|$ and that $\DeltaN$ is the set $\bx\in\Delta$, such that $Nm_s$ is an integer for all $s\in\{1\dots d\}$. For all $\bs\in\calS^N$, $\xs(\bs)\in\DeltaN$. Respectively, for all $\bx\in\DeltaN$, there exists $\bs\in\calS^N$ such that $\xs(\bs)=\bx$ (unique up to a permutation of its components).

    Our proof of the result uses the following lemma, which shows that the property holds for some particular values of $\bx$ and $\tilde{\bx}$:
    \begin{lem}
        \label{lem:intermediate_Lipschitz}
        The result of Lemma~\ref{lem:Lipschitz-lemma} holds for all $\bx,\tilde{\bx}$ for which there exists a $N>0$ such that $\bx,\tilde{\bx}\in\DeltaN$ and $\norme{\bx-\tilde{\bx}}\le 2/N$.
    \end{lem}

    Lemma~\ref{lem:intermediate_Lipschitz} implies that the result of Lemma~\ref{lem:Lipschitz-lemma} holds for all $N$ and for all $\bx,\tilde{\bx}\in\DeltaN$. Indeed, take any $\bx,\tilde{\bx}\in\DeltaN$, there exists an integer $k \in \mathbb{N}$ such that $\norme{\bx - \tilde{\bx}} = 2k/N$. We can then construct a sequence of $k+1$ state vectors $\bx = \bx_1$, $\bx_2$, $\dots$, $\bx_k$, $\bx_{k+1} = \tilde{\bx}$, such that $\norme{\bx_i - \bx_{i+1}} = 2/N$ for all $1 \le i \le k$. Hence
    \begin{align*}
      |\Vrel(\bx,t)-\Vrel(\tilde{\bx},t)|&  \le \sum_{i=1}^{k} |\Vrel(\bx_i,t)-\Vrel(\bx_{i+1},t)| \\
      & \le (1+(1-\gamma) + \dots + (1-\gamma)^{t-1})\Rmax \sum_{i=1}^{k} \norme{\bx_i - \bx_{i+1}}  \\
      & =  (1+(1-\gamma) + \dots + (1-\gamma)^{t-1})\Rmax\norme{\bx-\tilde{\bx}}
    \end{align*}
    By continuity of the value function and since $\cup_{N>0}\DeltaN$ is dense in $\Delta$, we can conclude that the result holds for all $\bx,\tilde{\bx}\in\Delta$.\qed

\end{proof}

Let us now return to the proof the Lemma \ref{lem:intermediate_Lipschitz}.

\begin{proof}{Proof of Lemma~\ref{lem:intermediate_Lipschitz}}
        We prove the result by induction on $t$. For $t=0$, the result is clearly true as $V_t(\bx)=0$ for all $\bx\in\Delta$. 
        
        We now assume that it holds for $t-1\ge0$ and look at what happens for $t$. Let $\bx,\tilde{\bx}\in\DeltaN$ be two population vectors such that $\norme{\bx-\tilde{\bx}}\le2/N$ and $\bx\ne\tilde{\bx}$. This implies that there exist two states $i\ne j\in\{1\dots d\}$ such that $x_i=\tilde{x}_i+1/N$ and $x_j=\tilde{x}_j-1/N$ (all others being equal), so that actually $\norme{\bx-\tilde{\bx}} = 2/N$. We denote by $\bs,\tilde{\bs}\in\calS^N$ two vectors such that $s_1=i$, $s'_1=j$, $s_n=s'_n$ for $n\ge2$, and $\xs(\bs)= \bx$, $\xs(\tilde{\bs})=\tilde{\bx}$.
        
        For a given vector of states $\bs\in\calS^N$ and a corresponding vector of actions $\ba=a_1\dots a_N$, with each $a_n$ representing a probability measure on $\calA$ as part of the randomized policy for arm $n$ in the relaxed problem, we define $Q(\bs,\ba)$ as the $Q$-value obtained by executing the action vector $\ba$ for the initial time step and subsequently adhering to the optimal policy:
        \begin{align}
            \label{eq:bellman-equation}
            Q(\bs,\ba) := \frac1N \sum_{n=1}^N r^{a_n}_{s_n} + \sum_{\bs'\in\calS^N} \Vrel(\xs(\bs'),t-1) \prod_{n=1}^{N} P^{a_n}_{s_n,s_n'}.
        \end{align}
        In the formulation above, $r^{a_n}_{s_n}$ and $P^{a_n}_{s_n,s_n'}$ should be interpreted as the \emph{expected} rewards and transition probabilities under the randomized action $a_n$ across $\calA$, that is $r^{a_n}_{s_n}=\sum_a r^{a}_{s_n}\proba{a_n=a}$ and $P^{a_n}_{s_n,s_n'}=\sum_a P^{a}_{s_n,s_n'}\proba{a_n=a}$. Consider $\ba^*$ as the optimal vector of actions for the state $\bs$, such that $\Vrel(\bx)=Q(\bs,\ba^*)$. Let $\ba'$ denote a vector of actions where $a'_1=0$ with probability one, and for all $n\in\{2\dots N\}$, $a'_n = a^*_n$. The action vector $\ba'$ represents an \emph{admissible} set of actions for $\tilde{\bs}$, given that action $0$ does not utilize any resources. However, it may not be optimal. Thus, it follows that $\Vrel(\tilde{\bx})\ge Q(\tilde{\bs},\ba')$, leading to
        \begin{align}
            \Vrel(\bx,t) &- \Vrel(\tilde{\bx},t) \le Q(\bs, \ba^*) - Q(\tilde{\bs},\ba')\nonumber\\
            &= \frac1N (r^{a^*_1}_{s_1}-r^{0}_{\tilde{s}_1}) + \sum_{\bs'\in\calS^N} \Vrel(\xs(\bs'),t-1) (\prod_{n=1}^{N}P^{a^*_n}_{s_n,s_n'} - \prod_{n=1}^{N}P^{a'_n}_{\tilde{s}_n,s_n'}) \nonumber\\
            &= \frac1N (r^{a^*_1}_{s_1}-r^{0}_{\tilde{s}_1}) + \sum_{\bs'\in\calS^N} \prod_{n=2}^{N}P^{a^*_n}_{s_n,s_n'}\Vrel(\xs(\bs'),t-1) (P^{a^*_1}_{s_1,s_1'}-P^0_{\tilde{s}_1,s'_1})  \nonumber\\
            &= \frac1N (r^{a^*_1}_{s_1}-r^{0}_{\tilde{s}_1}) + \sum_{s'_2, \dots, s'_N \in \calS^{N-1} } \prod_{n=2}^{N}P^{a^*_n}_{s_n,s_n'}\sum_{k\in\calS}\Vrel(\xs(k,s'_2\dots s'_N),t-1) (P^{a^*_1}_{s_1,k}-P^0_{\tilde{s}_1,k}),
            \label{eq:bounding_Lipschitz}
        \end{align}
        where we have used Equation~\eqref{eq:bellman-equation} for the first equality and the fact that $\bs$ and $\tilde{\bs}$, and $\ba$ and $\ba'$ are equal for all components except the first one. The last equality is a rewriting of the above equation where the special role of the first component has been emphasized.

        The first term of Equation~\eqref{eq:bounding_Lipschitz} is smaller than $2\Rmax/N$. For the second term, denote by $\gamma_{s_1 \tilde{s}_1, k} := \min ( P^{a^*_1}_{s_1, k}, P^{0}_{\tilde{s}_1, k})$. For any $s'_2, \dots, s'_N \in \calS^{N-1}$, we have 
        \begin{align}
            \sum_{k \in \calS} \Vrel&(\xs(k, s'_2\dots s'_N),t-1)(P^{a^*_1}_{s_1, k}-P^0_{\tilde{s}_1, k})\nonumber\\
            = & \sum_{k \in \calS} \Vrel(\xs(k, s'_2\dots s'_N),t-1)(P^{a^*_1}_{s_1, k} - \gamma_{s_1 \tilde{s}_1, k})
            -\sum_{k \in \calS} \Vrel(\xs(k, s'_2\dots s'_N),t-1)(P^0_{\tilde{s}_1, k} - \gamma_{s_1 \tilde{s}_1, k})\nonumber\\
            \le & (1 - \sum_{k \in \calS} \gamma_{s_1 \tilde{s}_1, k}) \left( \max_{\overline{k}} \Vrel(\xs(\overline{k}, s'_2\dots s'_N),t-1)
            - \min_{\underline{k}} \Vrel(\xs(\underline{k}, s'_2\dots s'_N),t-1) \right) \nonumber \\
            \le & (1 - \gamma) \left( \max_{\overline{k}} \Vrel(\xs(\overline{k}, s'_2\dots s'_N),t-1)
            - \min_{\underline{k}} \Vrel(\xs(\underline{k}, s'_2\dots s'_N),t-1) \right)
            \label{eq:bounding_Lipschitz2}
        \end{align}
        where the last step follows from our definition of $\gamma$. 
        
         Since $\xs(\overline{k}, s'_2\dots s'_N)$ and $\xs(\underline{k}, s'_2\dots s'_N)$ are at $\mathcal{L}^1$-distance $2/N$ when $\overline{k} \ne \underline{k}$, and $0$ otherwise, we can apply our induction hypothesis to Equation~\eqref{eq:bounding_Lipschitz2} and deduce from Equation~\eqref{eq:bounding_Lipschitz} that
        \begin{align*}
            \Vrel(\bx,t) &- \Vrel(\tilde{\bx},t) \le \frac{2\Rmax}{N} \left( 1+ (1-\gamma)(1+(1-\gamma)+ \dots + (1-\gamma)^{t-2}) \right)
        \end{align*}
        By exchanging the roles of $\bx$ and $\bx'$ and using that $\norme{\bx-\tilde{\bx}} = 2/N$, we finally obtain
        \begin{align*}
            |\Vrel(\bx,t) - \Vrel(\tilde{\bx},t)|& \le \frac{2\Rmax}{N}( 1+ (1-\gamma)(1+(1-\gamma) + \dots + (1-\gamma)^{t-2}))\\
            & = ( 1+ (1-\gamma) + \dots + (1-\gamma)^{t-1})\Rmax\norme{\bx-\tilde{\bx}},
        \end{align*}
        which completes the inductions step.
        \qed
    \end{proof}

We remark that to obtain a Lipschitz constant of $\Vrel(\cdot,t)$ that does not dependent of $t$, Lemma~\ref{lem:Lipschitz-lemma} has imposed the condition that $\gamma>0$. The constant $\gamma$ corresponds to the minimum probability that two arms starting in different states couple at the next time-step, when the action $0$ is applied to the first arm and any action is applied to the second arm. This condition can be weakened substantially. In particular, a less demanding sufficient condition is that there exists an integer $\tau \ge 1$ such that $\lambda^\tau>0$, where 
\begin{align*}
  \lambda^\tau := \min_{i \neq j, a_1\dots a_\tau} \sum_k \min( ( \underbrace{P^0\dots P^0}_{\tau \text{ times}} )_{ik}, (P^{a_1}\dots P^{a_\tau})_{jk}).
\end{align*}
The quantity $\lambda^\tau$ is the minimum probability that two arms starting in two different states couple after $\tau$ time-steps, where one arm applies the constant action $0$, while the other arm chooses arbitrary actions. Similar to Lemma~\ref{lem:Lipschitz-lemma}, we can prove the following result: 
\begin{lem}
  \label{lem:Lipschitz-lemma-tau}
  For all $\bx,\tilde{\bx}\in\Delta$, and all $\tau$ such that $\lambda^{\tau}>0$, the value function satisfies: 
  \begin{align*}
      |\Vrel(\bx,t)-\Vrel(\tilde{\bx},t)| \le \frac{\tau\Rmax}{\lambda^\tau}\norme{\bx-\tilde{\bx}}.
  \end{align*}
\end{lem}
\begin{proof}{Sketch of proof.}
  The proof of this results follows the same line as the proof of Lemma~\ref{lem:Lipschitz-lemma} and we only highlight the main steps:
  \begin{itemize}
    \item As in the proof Lemma~\ref{lem:Lipschitz-lemma}, the first step is to restrict our attention to points that are at distance $2/N$ in $\DeltaN$, which by using the representation $\bs$ implies that only one of the components differs between $\bs$ and $\tilde{\bs}$. 
    \item The main step of the proof of Lemma~\ref{lem:Lipschitz-lemma} is to write $\Vrel(\bx,t)$ as a function of $\Vrel(\bx,t-1)$ as in Equation~\eqref{eq:bellman-equation}, and uses Equations~\eqref{eq:bounding_Lipschitz} and \eqref{eq:bounding_Lipschitz2} to show that if $\norme{\bx-\tilde{\bx}}\le2/N$, then 
        \begin{equation*}
          |\Vrel(\bx,t)-\Vrel(\tilde{\bx},t)|\le \frac{2\Rmax}{N} + (1-\gamma)\max_{\bx',\tilde{\bx'}\text{ s.t. }\norme{\bx'-\tilde{\bx}'}\le2/N} |\Vrel(\bx',t-1)-\Vrel(\tilde{\bx}',t-1)|.
        \end{equation*}
    To prove Lemma~\ref{lem:Lipschitz-lemma-tau}, we can use a similar reasoning to show that if $\norme{\bx-\tilde{\bx}}\le2/N$, then
    \begin{align*}
      |\Vrel(\bx,t)-\Vrel(\tilde{\bx},t)|\le \frac{2\tau\Rmax}{N} + (1-\lambda^\tau)\max_{\bx',\tilde{\bx}'\text{ s.t. }\norme{\bx'-\tilde{\bx}'}\le2/N} |\Vrel(\bx',t-\tau)-\Vrel(\tilde{\bx}',t-\tau)|.
    \end{align*}
  \end{itemize}
  The result then follows from the same reasoning.   \qed
\end{proof}

\subsection{Non-degenerate problem and concentration on a trajectory}

The previous lemma can be extended to show that the stochastic system $\bXN(t)$ concentrates on a neighbourhood of $\bx^*(t)$.

\begin{lem}
  \label{lem:concentration_trajectory}
  Assume that the problem is non-degenerate, let $\by^*$ be the optimal solution to the LP computed at time-step $0$ and $\bXN(t)$ be the sequence of configuration vectors obtained when applying Algorithm~\ref{algo:LP-update-improved}. Then for all $\varepsilon>0$, there exists $C_1,C_2>0$ such that for all $N$:
  \begin{align}
    \label{eq:lemma_concentration_traj}
    \proba{\bXN(t)\in\calB(\bx^*(t),\varepsilon)} \ge 1-C_1 \cdot e^{-C_2 N}.
  \end{align}

  Moreover, if the LP has a unique solution starting from any initial point, then \eqref{eq:lemma_concentration_traj} also holds if $\bXN(t)$ is the output of  Algorithm~\ref{algo:LP-update}.
\end{lem}

\begin{proof}{Proof.}
  We first consider what happens when Algorithm~\ref{algo:LP-update-improved} is used. We proceed by induction on $t$. This is clearly true for $t=0$. Assume that this holds for some $t\ge0$. By Proposition~\ref{prop:optimal}, there exists $\varepsilon_t$ such that the control $\by(\bx)$ defined in Proposition~\ref{prop:optimal} is optimal for all $\bx\in\calB(\bx^*(t), \varepsilon_t)$. The induction hypothesis and the continuity of $\by(\bx)$ therefore imply that for all $\epsilon>0$, there exists $C_1,C_2>0$ such that $\norme{\bYN(t)-\by^*(t)}\le \epsilon$ with probability at least $1 - C_1 \cdot e^{-C_2 N}$, where $\bYN(t)$ is the control used for $\bXN(t)$.  Write
  \begin{align*}
   & \norme{\bXN(t+1) - \bx^*(t+1)} \le \norme{\bXN(t+1) - \phi(\bYN(t))} + \norme{\phi(\bYN(t))-\bx^*(t+1)}.
  \end{align*}
  Hence, by using Lemma~\ref{lem:Markovian-transition-analysis} and the union bound, for all $\epsilon>0$, there exists $C_1',C_2'>0$ such that $\norme{\bXN(t+1) - \bx^*(t+1)} \le \epsilon$ with probability at least $1-C_1' \cdot e^{C_2'N}$.

  To show \eqref{eq:lemma_concentration_traj}, the only remaining point is to see that if $s'$ is a state such that $x^*_{s'}(t+1)=0$, then so is $\XN_{s'}(t+1)$. By definition of the deterministic evolution $\bx^*(t+1)=\phi(\by^*(t))$ (see e.g. Equation~\eqref{eq:Markovian-transition}), we have:
  \begin{align*}
    x^*_{s'}(t+1) = \sum_{a,s} y^*_{s,a}(t)P^a_{s,s'}.
  \end{align*}
  Hence, $x^*_{s'}(t+1)$ equals $0$ if for all $s,a$, either $y^*_{s,a}(t)=0$ or $P^{a}_{s,s'}=0$. By construction of $\bYN(t)$ from $\by^*(t)$, $y^*_{s,a}(t)=0$ implies that $\YN_{s,a}(t)=0$, consequently $\XN_{s'}(t+1)=0$.

  To study what happens when Algorithm~\ref{algo:LP-update} is used, we remark that if the solution to the LP is unique, and $\bXN(t)$ is close enough to $\bx^*(t)$, then the new LP solution computed by Algorithm~\ref{algo:LP-update} is the same $\by(\bXN(t))$ used by Algorithm~\ref{algo:LP-update-improved}.  Hence, the proof for Algorithm~\ref{algo:LP-update-improved} also applies in this case.     \qed

\end{proof}

\subsection{Proof of the lower bounds}  \label{subsec:lowerbound}

Consider the following two-action restless bandit with two states $\calS=\{1,2\}$, with parameters $\bx(0)=(0.5,0.5)$, $T=2$, $A=\{0,1\}$, $\rr^0 = [0,0]$, $\rr^1 = [1,0]$, and where the transition matrices are
\begin{align*}
  \pp^0 = \begin{pmatrix}
    0.5 & \ 0.5 \\
    0.5 & \ 0.5
  \end{pmatrix}
  , \pp^1 =
  \begin{pmatrix}
    0.5 & \ 0.5 \\
    0.5 & \ 0.5
  \end{pmatrix}.
\end{align*}
We consider that $D(s,0) = 0$, $D(s,1) = 1$ for any state $s$, and we distinguish the resource constraints $b=0.3$ and $b=0.5$.

For any resource constraint, the solution to the LP is to choose action $1$ for as many arms as possible. This gives a reward $2b$. If $b=0.3$, the problem is non-degenerate whereas if $b=0.5$ the problem is degenerate.

For the stochastic system with $N$ components, this gives a reward $\frac1N\floor{Nb}$ at time-step $0$ and a reward $\min(b, \frac1N \floor{N \XN_1(1)})$ at time-step $1$. Since $\XN_1(1)$ follows a binomial distribution of parameter $(N,0.5)$, the total reward of the LP-update policy is equal to
\begin{align*}
  \frac1N\floor{Nb}+b + \expect{\min(\frac1N \floor{N \XN_1(1)}-b, 0)}.
\end{align*}
By the central limit theorem, $\lim_{N\to\infty}\sqrt{N}\expect{\min(\XN_1(1)-0.5, 0)}=\sqrt{2}{\pi}>0$. This provides a counter-example for the lower bound of Theorem~\ref{thm:LP-update-general}.

If $b=0.3$ and $N$ is not a multiple of $10$, the problem does not admit a perfect rounding, and we have $2b - (\frac1N\floor{Nb}+b + \expect{\min(\frac1N \floor{N \XN_1(1)}-b, 0)}) \ge b-\frac1N\floor{Nb}\ge 0.1/N$. This provides a counter-example for the lower-bound of Theorem~\ref{thm:LP-update-rate-nondegenerate}.

If $b=0.3$ and $N$ is a multiple of $10$, then the problem admits a perfect rounding. In this case, classical anti-concentration arguments (\citet{matouvsek2001probabilistic}) show that
\begin{align*}
  \proba{\XN_1(1)\le 0.2} \ge \frac1{15}\exp(-16N(0.6-0.2)^2).
\end{align*}
This shows that
\begin{equation*}
\expect{\min(\frac1N \floor{N \XN_1(1)}-0.3, 0)} \ge -0.1 \proba{\XN_1(1)\le 0.2} = -\frac1{150}\exp(-2.56N).
\end{equation*}
It provides a counter-example for the lower-bound of Theorem~\ref{thm:LP-update-rate-expo}.  \qed

\section{\bfseries\scshape{The occupation measure policy}} \label{subsec:occupation-measure-policy}

To provide a benchmark for evaluation of the performance of the LP-update policy, we describe in this section a policy that we call the \emph{occupation measure policy}.  This policy is a straightforward generalization of the randomized activation control policy in
\citet{ZayasCabn2017AnAO} and the occupancy-measured-reward index policy in \citet{xiong2021reinforcement} to the multi-constraint multi-action case.  The occupation measure policy is proven to be asymptotically optimal for the multi-action single-constraint multi-armed bandits in \citet{xiong2021reinforcement}. Adapting their proof to the multi-constraint case is direct. 

The occupation measure policy is a one-pass policy that solves the linear program only once at the very beginning, and constructs the probability vectors $\mu^*_{s,a}(t)$ defined in \eqref{eq:occupation-measure} from the solution. At each decision epoch $t$, the budget left is initialized as $\bB := N \cdot \bb$, the action on each arm $n$ is initialized as the passive $0$. It then samples a new action $a_n$ from the distribution  $(\mu^*_{s_n,a}(t))_{a \in A(s_n)}$ on each arm $n$. If choosing action $a_n$ instead of action $0$ on arm $n$ does not violate any of the budget constraints, then we apply action $a_n$ on arm $n$, and we decrease the budget left $\bB$; otherwise we keep action $0$ on arm $n$ and continue to sample an action on the next arm. The detailed implementation is given  in Algorithm \ref{algo:occupation_measure}.

\begin{algorithm}[t]
\SetAlgoLined
\SetKwInput{KwInput}{Input}
\KwInput{Time horizon $T$ and initial configuration vector $\bx(0)$.}
  Solve the linear program \eqref{eq:relaxed_problem_general} with time horizon $T$ and initial configuration vector $\bx(0)$, obtain an optimal solution $\by^*$ and the corresponding $\bx^*$ \;
 \For{$t = 0,1,2,\dots,T-1$}{
  Compute from the LP solution the occupation measure \begin{equation} \label{eq:occupation-measure}
            \mu^*_{s,a}(t) := \begin{cases}
                                \frac{y^*_{s,a}(t)}{x^*_s(t)}, & \mbox{if } x^*_s(t) > 0 \\
                                \mathbf{1}_{\{ a=0\} } , & \mbox{otherwise}.
                              \end{cases}
          \end{equation}
  Observe the current states of the $N$ sub-MDP's $\bs = (s_1,s_2, \dots, s_N)$. Initialize $\bB := N \cdot \bb$, and actions on the $N$ sub-MDP's as $A(\bs) = (0,0,\dots,0)$ \;
  \For{$n=1,2,\dots,N$}{ Sample an action $a_n$ according to the probability vector $(\mu^*_{s_n,a}(t))_{a \in A(s_n)}$ \;
      \If{$\bB - D(s_n,a_n) \ge \mathbf{0}$ (component-wise)}{
          $A(s_n) := a_n$ \;
          $\bB := \bB - D(s_n,a_n)$ \;

      }
  }
  Apply the actions $A(\bs)$ to the $N$ sub-MDP's \;
  }
 \caption{Occupation measure policy for weakly-coupled MDPs.}
 \label{algo:occupation_measure}
\end{algorithm}

\section{\bfseries\scshape{Comparison with existing notions of non-degeneracy}} \label{subsec:alternative-formulation-of-non-degeneracy}

The non-degenerate property as defined in Definition \ref{def:non-degeneracy} uses a rank condition from the optimal LP solution $\by^*$. In this section, we provide a detailed comparison of our notion to the existing notions presented \cite{zhang2021restless,gast2023linear,brown2023fluid,zhang2022near1}. To facilitate the comparison, we reformulate our definition. 

Let $\by^*$ be an optimal solution to the LP (at time $0$) and let us fix a time $t$.  For each state $s\in\calS$, we define by $A_s^+ := |\{a \in \calA \mid y^*_{s,a}(t) > 0 \}|$ the number of different actions that are chosen in state $s$ at time $t$. Similar to \citet{brown2023fluid}, we define:
\begin{align*}
  \calS^{+} & := \{ s \in \calS \mid {A_s^+} = 1  \}; \\
  \calS^{0} & := \{ s \in \calS \mid {A_s^+} > 1  \}; \\
  \calS^{-} & := \{ s \in \calS \mid {A_s^+} = 0  \}.
\end{align*}
We see that $\calS^{+}$, $\calS^{0}$, $\calS^{-}$ form a partition of $\calS$. Their interpretations are that (for the LP), a state in $\calS^{+}$ has a unique optimal action, a state in $\calS^{0}$ has at least two equally optimal actions, and a state in $\calS^{0}$ is not reached at time $t$.

As $\calU^*(t)$ is the set of indices $(s,a)$ for which $y^*_{s,a}(t) = 0$, we have: 
\begin{align*}
  \abs{\calU^*} &= \sum_{s \in \calS} \abs{\calA} - {A_s^+}\\
  &= \sum_{s \in \calS^{-}} {A_s} + \sum_{s \in \calS^{+}} {A_s} - \abs{\calS^{+}} + \sum_{s \in \calS^{0}} \abs{\calA} - \sum_{s \in \calS^{0}} {A_s^+}  \\
  & = \abs{\calU} - \abs{\calS^{+}} - \sum_{s \in \calS^{0}} {A_s^+},
\end{align*}
while $\abs{\calS^*} = \abs{\calS^{0}} + \abs{\calS^{+}}$. 

As $\abs{\calJ^*}$ is the number of saturated constraints, we have
\begin{equation*}
  |\calJ^*| + |\calS^*|+|\calU^*| = \abs{\calU} + |\calJ^*| + \abs{\calS^{0}} - \sum_{s \in \calS^{0}} {A_s^+}.
\end{equation*}

The matrix $C^*$ has $|\calJ^*| + |\calS^*|+|\calU^*|$ rows and $|\calU|$ columns. To simplify the discussion, in this section, we assume that:
\begin{equation}\label{eq:additoinal-assumption}
  \mbox{All rows of $C^*(t)$ are linearly independent.}
\end{equation}
So that the rank of $C^*(t)$ boils down to verifying an inequality among dimension. Indeed, under this assumption \eqref{eq:additoinal-assumption}, in order for $C^*$ to have rank $|\calJ^*| + |\calS^*|+|\calU^*|$ as in Definition \ref{def:non-degeneracy}, it is necessary and sufficient that 
\begin{equation}\label{eq:reformulation-of-non-degeneracy}
  |\calJ^*| + \abs{\calS^{0}} \le \sum_{s \in \calS^{0}} {A_s^+}.
\end{equation}
 
We are now ready to compare the formulation \eqref{eq:reformulation-of-non-degeneracy} with various existing notions of non-degeneracy in the literature:
\begin{enumerate}
  \item In \citet{zhang2021restless} and \citet{gast2023linear} that study two-action bandits with a single \emph{equality} resource constraint $b = \alpha$. Their notion of non-degeneracy using our notation is to say that $\abs{\calS^{0}} \ge 1$. Since in their model the assumption \eqref{eq:additoinal-assumption} always holds, and $\abs{\calJ^*}$ is always $1$. In order for $\eqref{eq:reformulation-of-non-degeneracy}$ to hold, it is necessary and sufficient that $\abs{\calS^{0}} \ge 1$, which shows that the two definitions of non-degeneracy are indeed equivalent. 
  \item In \citet{brown2023fluid}, the authors define another notion of non-degeneracy in Definition 5.2 for multi-action single-constraint bandits, which by using the notations in this section can be stated as follows (assuming \eqref{eq:additoinal-assumption}): the problem is non-degenerate if at least one of the conditions hold: (i) the single resource constraint is not saturated; (ii) $\abs{\calS^{0}} \ge 1$. This is again equivalent to \eqref{eq:reformulation-of-non-degeneracy}. Note that additional care is taken in Definition~5.2 of \citet{brown2023fluid}, in order to account for the non-satisfaction of \eqref{eq:additoinal-assumption}. 
  \item In \citet{zhang2022near1}, another definition of non-degeneracy for multi-action multi-constraint bandits is given in Section 3.3. According to this definition, a multi-action multi-constraint bandit is non-degenerate if $\abs{\calJ^*} \le 1$ and $\abs{\calS^{0}} \ge 1$. Clearly this implies \eqref{eq:reformulation-of-non-degeneracy}, but under this general model, there are other possible problems that satisfy \eqref{eq:reformulation-of-non-degeneracy} while not covered by the non-degenerate definition of \citet{zhang2022near1}. 
\end{enumerate}


\section{\bfseries\scshape{Improving the rounding method}} \label{subsec:rounding-problem-discussion}

When the LP is solved with a given initial condition $\bXN$, it yields a first-step control $\bY:=\by(0)\in\calY(\bXN)$. However, this solution may not be directly transferable to the original problem because the quantities $N \cdot Y_{s,a}$ are not necessarily integers. Throughout our paper, we address this issue by employing a simple truncation method to obtain a vector $\bYN$ that is defined as follows:
\begin{equation} \label{eq:truncation}
  \bYN_{s,a} = \left\{
  \begin{array}{ll}
    N^{-1}\floor{NY_{s,a}}  &\mbox{ if }  a \not=0\\
    \XN_{s} - \sum_{a\not = 0} \tilde{Y}_{s,a} & \mbox{ if }  a =0.
  \end{array}
\right.
\end{equation}
This approach results in an admissible solution $\bYN\in\calYN(\bXN)$, as all elements in $D(s,a)$ and $\bb$ are non-negative, and notably, $D(s,0) = 0$. However, this truncation solution introduces an $\calO(1/N)$ sub-optimality gap.

In this section, we delve into potential strategies to enhance performance and reduce this sub-optimality gap. By examining alternative methods to handle the non-integer nature of the solution, we aim to improve the applicability and efficiency of our approach, ensuring closer alignment with the optimal solution of the original problem.

\subsection{Improvement by using MILP}   \label{subsec:sophisticated-rounding}

As previously mentioned, the rounding process we employ in \eqref{eq:truncation} addresses the issue that the first step $\by(0)$ from the LP might not yield integer values for $N \cdot y_{s,a}(0)$. To circumvent this, another more refined approach is to explicitly incorporate integer constraints into our LP, leading to the following mixed-integer optimization problem (MILP):

\begin{maxi!}|s|{\by\ge0}{\sum_{t=0}^{T-1} \sum_{(s,a) \in \calU} r_s^a y_{s,a} (t)}{\label{eq:relaxed_problem_general_integer}}{}
  \addConstraint{\sum_{a \in \calA} y_{s,a}(0) = x_s}{}{\forall s \label{eq:init_general_integer}}
  \addConstraint{\sum_{a \in \calA} y_{s,a}(t+1) = \big( \phi(\by(t)) \big)_s \label{eq:markov_relaxed_general_integer}\qquad}{}{\forall s,t}
  \addConstraint{D \by(t) \le \bb \label{eq:relaxed-constraints_integer}}{}{\forall t},
  \addConstraint{N \cdot y_{s,a}(0)\in\N \label{eq:integer}}{}{\forall s,a},
\end{maxi!}
It is important to note that since the solution primarily utilizes control at the first time-step, we apply the integrity constraint solely at $t=0$.

We introduce $\UpdateMILP$ to denote the expected performance of the policy obtained by implementing Algorithm~\ref{algo:LP-update}, wherein the solution of the mixed-integer linear program \eqref{eq:relaxed_problem_general_integer} is utilized instead of \eqref{eq:relaxed_problem_general}. While this modified problem is more constrained than the original LP, it still represents a relaxation of the initial problem. Therefore, for any admissible policy $\pi\in{\mathrm{LP-update},\mathrm{MILP-update}}$, the following holds:
\begin{align*}
  V_\pi^{(N)}(\bx,T) \le \Vrel(\bx,T)
\end{align*}
In Figure~\ref{fig:MIP}, we present a performance comparison between the LP-update and MILP-update policies. As anticipated, the MILP-update demonstrates superior performance compared to the original LP-update, particularly for smaller values of $N$. 

\begin{figure}[ht]
  \centering
  \includegraphics[width=0.6\linewidth]{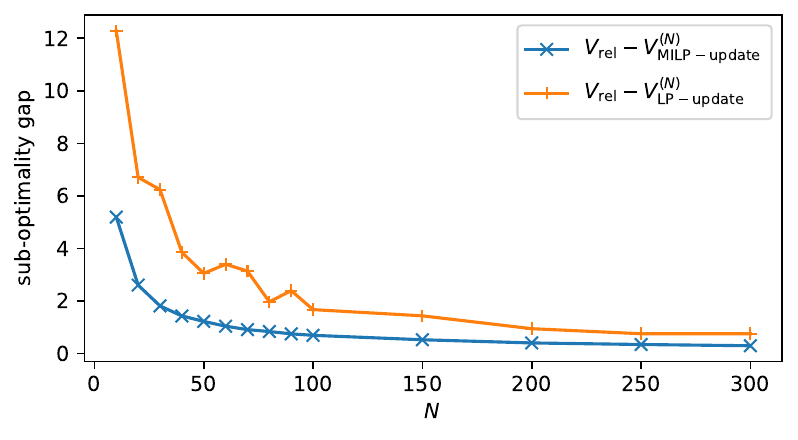}
  \caption{Sub-optimality gap of the LP-update policy compared to the MILP-update policy.}
  \label{fig:MIP}
\end{figure}

From a computational perspective, addressing \eqref{eq:relaxed_problem_general_integer} is more demanding in terms of resources compared to solving the original LP problem. Our experimental findings revealed that solving the MILP required no more than twice the time needed for the original LP for this specific model having $10$ states. Notably, the duration for solving the MILP appeared to be unaffected by the size of $N$, suggesting that increases in $N$ do not proportionally complicate the MILP solution process. This consistency in computational time, regardless of the scale of $N$, offers an interesting topic for further investigation, especially to determine if this trend is consistent across various scenarios and applications. It is also important to observe that the zigzag effect in the LP-update policy, resulting from the cutoff rounding, is absent in the MILP-update policy. This suggests that the MILP-update policy might be a more suitable and natural choice for implementation.

\subsection{Perfect rounding}



The constructed $\bYN$, even using MILP, might still be far from $\bY$ for small values of $N$.  To obtain a $\bYN$ closer to $\bY$, the approach developed in \citet{gast2023linear} for two-action bandits is to use the so-called \emph{randomized rounding}, that consists of using a random vector $\bYN$ such that $\expect{\bYN}=\bY$ and $\bYN\in\calYN(\bXN)$ almost surely. If such a random vector exists, we say that the problem admits a \emph{perfect rounding}.
\begin{defi}[Perfect rounding]\label{def:perfect-rounding} \
  An LP problem admits a perfect rounding for an integer $N$, if for all $\bx\in\DeltaN$ and for all $\by\in\calY(\bx)$, there exists a random vector $\bYN$ such that $\bYN\in\calYN(\bx)$ almost surely and $\expect{\bYN} = \by$.
\end{defi}
The next theorem shows that, if there exists a perfect rounding, we obtain the much faster exponential convergence rate by avoiding the rounding errors. The result is summarized in the following convergence theorem:
\begin{thm} \label{thm:LP-update-rate-expo}  
  Denote by $\Update$ the value of the LP-update policy defined in Algorithm \ref{algo:LP-update} or \ref{algo:LP-update-improved}, and by $\rel$ the value of the linear program \eqref{eq:relaxed_problem_general}. Then
  \begin{itemize}
    \item If the weakly-coupled MDP with statistically identical arms is non-degenerate, then there exist constants $C_1,C_2>0$ such that for any $N$ for which the problem admits a perfect rounding:
          \begin{equation*}
            \abs{\Update - \rel} \le C_1 \cdot e^{-C_2 N}.
          \end{equation*}
    \item There exists a non-degenerate weakly-coupled MDP with statistically identical arms that admits a perfect rounding for an infinite number of $N$, and two constants $C_1',C_2'$ such that for all such $N$:
    \begin{equation*}
      \abs{\Update - \rel} \ge C_1' \cdot e^{-C_2' N}.
    \end{equation*}
  \end{itemize}
\end{thm}

\begin{proof}{Proof.}
  The proof of the first item is identical to the proof of Theorem~\ref{thm:LP-update-rate-nondegenerate} until the last equation, for which if we have a perfect rounding, then $\expect{\bYN(t)-\bY(t)}=0$. The proof of the lower bound in the second item is provided in Section \ref{subsec:lowerbound}. \qed
\end{proof}

\subsubsection*{Comparison with existing literature}

The authors of \cite{gast2023linear} study the single constraint case where $\calA=\{0,1\}$, $D(s,a)=a$ and $b=\alpha$. They show that a perfect rounding exists for all $N$ such that $\alpha N$ is an integer. Another example is  when  the constraints  \eqref{eq:positive4}-\eqref{eq:relaxed-constraints-1}-\eqref{eq:init_general-1}
form a totally unimodular matrix (see for instance \citet{DBLP:books/daglib/p/HoffmanK10}) for all $t$: in that case, the solution to the LP satisfies $\by\in\calYN(\bXN)$ automatically and no rounding is needed. There are, however, many situations for which perfect rounding is impossible. In such cases, a method to improve the simple truncation \eqref{eq:truncation} is to find $\bYN\in\calYN(\bXN)$ that minimizes the distance $\snorme{\bYN-\bY}$. This can be computed by solving an integer linear program.

The researchers in \cite{ZayasCabn2017AnAO,xiong2021reinforcement} implement a randomized approach to derive an admissible $\bYN$ from $\by^*$. This method is detailed in our description of their algorithm in Algorithm~\ref{algo:occupation_measure}. Initially, their algorithm sets all arms to action $0$. Subsequently, it iteratively processes each arm. For an arm in state $s$, the algorithm probabilistically selects a new action $a$ based on the ratio $y^*_{s,a}(t)/x^*_s(t)$. This action $a$ is then assigned to the arm, provided it does not breach any budget constraints. The central limit theorem underpins the rationale for this approach, ensuring that $\bYN$ approximates $\by^* + \calO(1/\sqrt{N})$. However, it does not extend to achieving a tighter approximation of $\by^* + \calO(1/N)$. Consequently, even in the context of non-degenerate problems, the randomized process employed by these algorithms results in a sub-optimality gap of $\calO(1/\sqrt{N})$ rather than the refined $\calO(1/N)$.

\end{APPENDICES}

\end{document}